# Minimax Optimal Estimation in Partially Linear Additive Models under High Dimension

Zhuqing Yu*, Michael Levine† and Guang Cheng‡
January 7, 2018

*Department of Statistics*
*Purdue University*
*250 N. University Street*
*West Lafayette, IN 47906*
*Email: zhuqing.yu.stat@gmail.com*
*mlevins@purdue.edu*
*chengg@stat.purdue.edu*

In this paper, we derive minimax rates for estimating both parametric and nonparametric components in partially linear additive models with high dimensional sparse vectors and smooth functional components. The minimax lower bound for Euclidean components is the typical sparse estimation rate that is independent of nonparametric smoothness indices. However, the minimax lower bound for each component function exhibits an interplay between the dimensionality and sparsity of the parametric component and the smoothness of the relevant nonparametric component. Indeed, the minimax risk for smooth nonparametric estimation can be slowed down to the sparse estimation rate whenever the smoothness of the nonparametric component or dimensionality of the parametric component is sufficiently large. In the above setting, we demonstrate that penalized least square estimators can nearly achieve minimax lower bounds.

*AMS 2000 subject classifications:* Primary 62C20, 62G05, 62F10; secondary 62G86, 62J12 .
*Keywords:* High dimension; minimax optimal; partial linear additive model; semiparametric.

## 1. Introduction

In this paper, we consider high dimensional partially linear additive models:

$$Y = X^T\beta_0 + \sum_{j=1}^{J} f_j(Z_j) + \varepsilon, \tag{1.1}$$

where the Euclidean vector $\beta_0 \in \mathbb{R}^p$ is sparse with $p > n$ and $f_j : \mathbb{R} \mapsto \mathbb{R}$ are nonparametric functions with possibly different smoothness. Assume $J$ is fixed while sparsity and smoothness parameters are known. Under this setting, minimax risks of estimation for both components are derived. As a side note, we mention that the choice of model structure, i.e., which covariate is linear or nonlinear, can be determined by the method developed in Zhang et al. (2011).

Without loss of generality, we assume $J = 2$ in this paper:

$$Y = X^T\beta_0 + f_0(Z) + g_0(U) + \varepsilon, \tag{1.2}$$

where $\beta_0 \in \mathbb{R}^p$ has at most $s_0$ non-zero elements, and $f_0$ and $g_0$ belong to the $\alpha$-th and $\gamma$-th order Sobolev balls, respectively. The $\alpha$-th order Sobolev ball over $[0,1]$, denoted as $W^{\alpha,2}(L_1)$, is defined as

*PhD Student. Supported by NSF DMS-1151692, DMS-1712919 and Purdue Research Foundation. Currently employed at AbbVie Inc., North Chicago, IL

†Associate Professor. Research Sponsored by NSF-DMS 1208994.

‡Professor. Research Sponsored by NSF CAREER Award DMS-1151692, DMS-1712919, Simons Fellowship in Mathematics and Office of Naval Research (ONR N00014-15-1-2331).

$\{f : [0,1] \to \mathbb{R} | J_\alpha^2(f) \le L_1^2\}$ for a constant $L_1 > 0$, where $J_\alpha^2(f) = \int_0^1 (f^{(\alpha)}(z))^2\, dz$ with $f^{(\alpha)}$ being the $\alpha$-th derivative of $f$. When the dimension of $\beta_0$ is fixed or slowly increasing ($p < n$), the above model has been extensively studied in the semiparametric literature, e.g., Härdle et al. (2000); Xie and Huang (2009); Cheng et al. (2015), while the high dimensional extension with $p > n$ has been further considered in Müller and van de Geer (2015); Ma and Huang (2016); Zhu (2017). Despite these literature, the minimax rates in estimating $\beta_0$, $f_0$ and $g_0$ remain unclear as far as we are aware.

First, we establish the minimax lower bound for estimating $\beta_0$ as

$$R_{\beta_0}(s_0, \alpha, \gamma) \gtrsim \frac{s_0}{n} \log\left(\frac{p}{s_0}\right), \tag{1.3}$$

up to a universal constant, based on iid observations $\{Y_i, X_i, Z_i, U_i\}_{i=1}^n$. It is worth noting that the lower bound does not depend on nonparametric smoothness indices, say $\alpha$ and $\gamma$, and coincides with the classical sparse estimation rate in the high dimensional linear models (Ye and Zhang, 2010; Raskutti et al., 2011; Verzelen, 2012). This result is similar in spirit to the classical low dimensional result where the Euclidean part can be estimated at $\sqrt{n}$-rate even in the presence of nuisance functions with slower rates; see Bickel et al. (1993).

A somewhat surprising result is that the lower bounds for estimating $f_0$ and $g_0$ turn out to be affected by the existence of $\beta_0$:

$$R_{f_0}(s_0, \alpha, \gamma) \gtrsim \max\left(n^{-2\alpha/(2\alpha+1)}, \frac{s_0}{n} \log\left(\frac{p}{s_0}\right)\right), \tag{1.4}$$

and

$$R_{g_0}(s_0, \alpha, \gamma) \gtrsim \max\left(n^{-2\gamma/(2\gamma+1)}, \frac{s_0}{n} \log\left(\frac{p}{s_0}\right)\right). \tag{1.5}$$

Such one-way interaction can be intuitively explained by the orthogonal decomposition (2.6). An interesting consequence of (1.4) and (1.5) is that the best possible estimation of $f_0$ and $g_0$ could be slowed down to the well known sparse estimation rate. To demonstrate this rate-switching phenomenon, we plot a two regime dichotomy in Figure 1: (i) in the *sparse* regime where $f_0$ is sufficiently smooth or $p$ is sufficiently high, the minimax risk lower bound becomes $s_0 \log(p/s_0)/n$; (ii) in the *smooth* regime where $f_0$ is very rough or $p$ is low, the lower bound becomes the classical nonparametric rate $n^{-2\alpha/(2\alpha+1)}$ (Pinsker, 1980; Stone, 1985). Note that a similar phase transition phenomenon occurs in high dimensional additive nonparametric models but due to very different reasons; see Koltchinskii and Yuan (2010); Raskutti et al. (2012); Suzuki and Sugiyama (2013); Yuan and Zhou (2016). We also note that the lower bound of estimating $f_0$ or $g_0$ does not depend on the smoothness of the other nonparametric component. This result essentially generalizes Horowitz et al. (2006) who showed that, in an additive nonparametric regression model, each component can be estimated (up to the first order asymptotics) as well as if all the rest were known.

In contrast with the literature on sparse parametric or nonparametric estimation such as Koltchinskii and Yuan (2010); Ye and Zhang (2010); Raskutti et al. (2011, 2012); Suzuki and Sugiyama (2013); Yuan and Zhou (2016), we are not interested in estimating the conditional mean function $E(Y|X, Z_1, \ldots, Z_J)$ as a whole, but rather separate minimax risk for each model component: $\beta_0, f_0, g_0$. Note that our results are not directly implied by the above papers where additive components are always assumed to share the same linear or nonlinear structure with the same smoothness.

In the end, we demonstrate that the penalized least square estimate for $(\beta_0, f_0, g_0)$ can almost achieve the lower bounds established above. To obtain such estimation rates, we develop a series of oracle inequalities that give more and more refined estimation errors for each model component in the order of $g, f$ and $\beta$ (under the assumption that $f$ is smoother than $g$), and then derive the risk upper bounds by strengthening these oracle inequalities to their moment versions.

**Notations**. For any vector $v \in \mathbb{R}^n$, we write its $\ell_1$, Euclidean and $\ell_\infty$ norm as $\|v\|_1 = \sum_{i=1}^n |v_i|$, $\|v\| = \sqrt{\sum_{i=1}^n v_i^2}$ and $\|v\|_\infty = \max_{1\le i\le n} |v_i|$, respectively, and also $\|v\|_n^2 := v^T v/n$. With a bit abuse of notation, we define for any function $f : \mathcal{Z} \mapsto \mathbb{R}$ that $\|f\| = \sqrt{\mathbb{E} f^2(Z)}$, $\|f\|_\infty = \sup_{z\in[0,1]} |f(z)|$ and

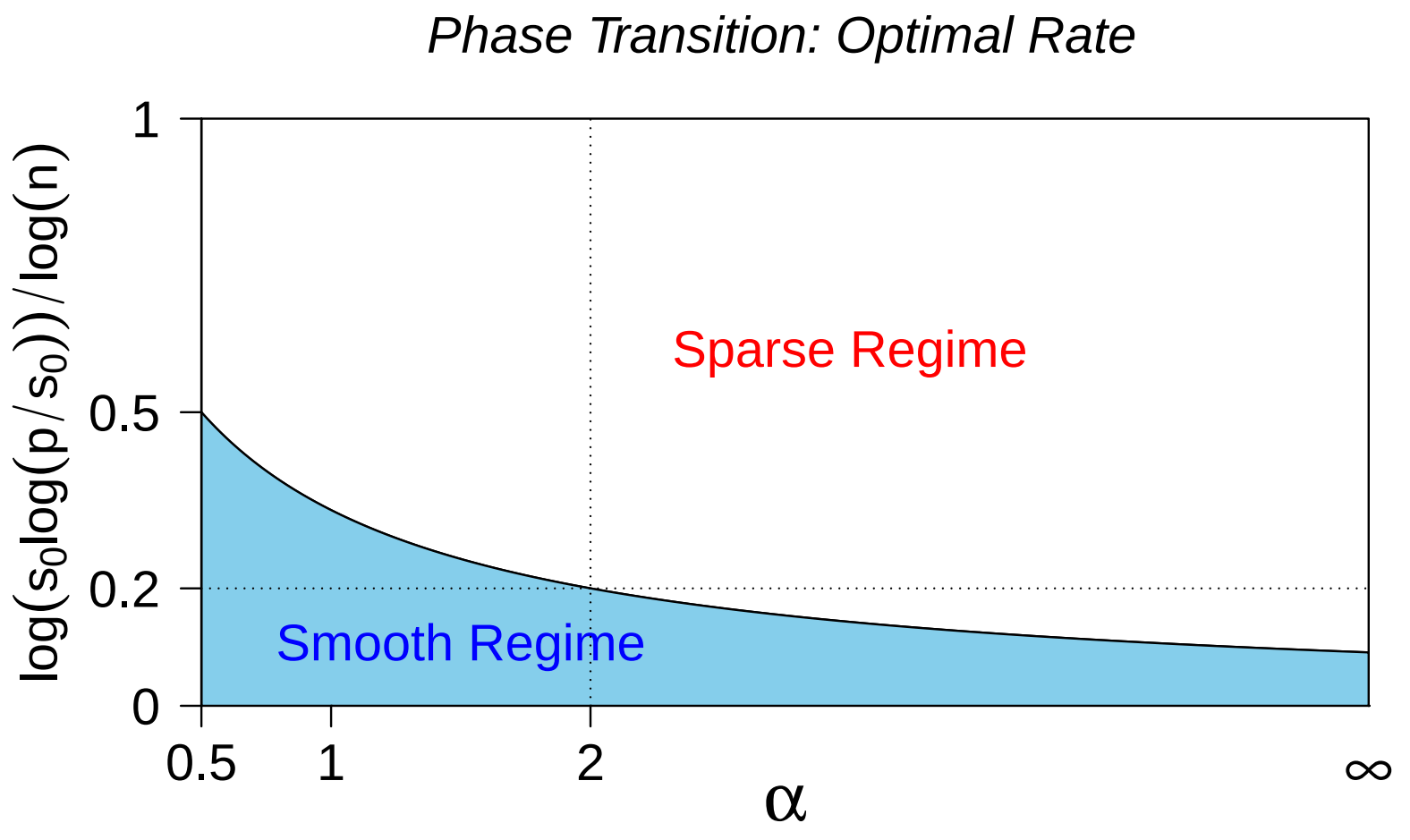

**Figure 1**. The minimax lower bound is $n^{-2\alpha/(2\alpha+1)}$ when $\alpha, p, s_0$ and $n$ fall into smooth regime. Otherwise, the minimax lower bound is $s_0 \log(p/s_0)/n$ in the sparse regime.

$\|f\|_n^2 = \sum_{i=1}^n f^2(Z_i)/n$. Let $S_0$ be the set of all non-zero components of $\beta_0$ and $s_0 = |S_0|$. Define $\beta_{S_0}$ such that $(\beta_{S_0})_j = \beta_j \mathbf{1}\{\beta_{0j} \neq 0\}$ and $\beta_{S_0^c} = \beta - \beta_{S_0}$, for any $\beta \in \mathbb{R}^p$. Thus, $\|\beta\|_1 = \|\beta_{S_0}\|_1 + \|\beta_{S_0^c}\|_1$. For any $x \in \mathbb{R}$, $\lceil x \rceil$ is the smallest integer that is strictly greater than $x$. For real sequences $a_n, b_n$, if $a_n \lesssim b_n$ ($a_n \gtrsim b_n$), then $\limsup a_n/b_n \leq C$ ($c \leq \limsup a_n/b_n$), for some constant $C$ (constant $c$). If $a_n \asymp b_n$, then $c \leq \liminf a_n/b_n \leq \limsup a_n/b_n \leq C$ for some constant $c, C$. Also, we write $a_n = O(b_n)$ if $|a_n| \leq C|b_n|$ for some constant $C > 0$. In the sequel, $c, c', C, C', \ldots$ denote a generic constant which may differ at each appearance.

# 2. Main Results

## 2.1. Minimax Lower Bounds

In this section, we assume $X$ is a mean zero Gaussian vector with variance matrix $\Sigma$, and the errors $\{\varepsilon_i\}_{i=1}^n$ are i.i.d. standard Gaussian random variables independent of $\{X_i, Z_i, U_i\}_{i=1}^n$. For simplicity, we standardize $X$ such that the diagonal of $\Sigma$ consist of 1's. Under this setting, we establish separate lower bounds on the minimax risk of estimating $\beta_0$, $f_0$ and $g_0$. For identifiability purpose, we assume $\mathbb{E} f_0(Z) = 0$.

We are ready to define the risk for estimating $\beta_0$ as

$$R_{\beta_0}(s_0, \alpha, \gamma, \Sigma) := \inf_{\widehat{\beta}} \sup_{\substack{\beta_0 \in B[s_0,p], f_0 \in W^{\alpha,2}(L_1), \\ g_0 \in W^{\gamma,2}(L_2)}} \mathbb{E}[\|\beta_0 - \widehat{\beta}\|^2], \tag{2.1}$$

where $B[s_0, p]$ and $\mathcal{S}_p$ denote a set of $p$-dimensional vectors with at most $s_0$ non-zero coordinates and a set of $p \times p$ covariance matrices with 1's on the diagonal, respectively. Since the supremum of minimax risks with respect to all covariance matrices $\Sigma$ is $+\infty$, it only makes sense to consider the infimum of minimax risks with respect to random designs:

$$R_{\beta_0}(s_0, \alpha, \gamma) := \inf_{\Sigma \in \mathcal{S}_p} R_{\beta_0}(s_0, \alpha, \gamma, \Sigma),$$

as indicated by Verzelen (2012). Similarly, we define the risk of estimating $f_0$ as

$$R_{f_0}(s_0, \alpha, \gamma) := \inf_{\Sigma \in \mathcal{S}_p} R_{f_0}(s_0, \alpha, \gamma, \Sigma),$$

where

$$R_{f_0}(s_0, \alpha, \gamma, \Sigma) = \inf_{\widehat{f}} \sup_{\substack{\beta_0 \in B[s_0,p], f_0 \in W^{\alpha,2}(L_1), \\ g_0 \in W^{\gamma,2}(L_2)}} \mathbb{E} \int_0^1 |\widehat{f}(z) - f_0(z)|^2 \, dz.$$

$R_{g_0}(s_0, \alpha, \gamma)$ is defined similarly.

Our main result in this paper is on the minimax lower bound presented below. We start with a version of the Fano's Lemma, i.e., Corollary 2.19 in Massart (2007), to be used in the proof. Suppose that $s = (\beta, f)' \in S$ where $S = \mathbb{R}^p \times \mathcal{F}$. The induced probability measure is written as $P_s$. A finite subset of $\mathcal{F}$ is denoted as $\mathcal{C}_1$ and a finite subset of $\mathbb{R}^p$ is denoted as $\mathcal{C}_2$. Their Cartesian product is denoted as $\mathcal{C}$ with the obvious cardinality $|\mathcal{C}| = |\mathcal{C}_1||\mathcal{C}_2|$.

**Lemma 2.1.** *We consider a set of statistical models $\{P_s, s \in S\}$ where $(S, d)$ is a pseudometric space. Let $\kappa$ be the absolute constant suggested in Corollary 2.18 of Massart (2007). Choose an arbitrary estimator $\widehat{s} = (\widehat{\beta}, \widehat{f})$ of $s$ and a finite subset $\mathcal{C} = \mathcal{C}_1 \times \mathcal{C}_2$ of $S$, such that $\max_{s,t \in \mathcal{C}} K(P_s, P_t) \leq \kappa \log |\mathcal{C}|$. Then, setting $\delta = \min_{\substack{s,t \in \mathcal{C} \\ s \neq t}} d(s, t)$, we have for any $r \geq 1$*

$$\sup_{s \in \mathcal{C}} \mathbb{E}\left[d^r(s, \widehat{s})\right] \geq 2^{1-r} \delta^r (1 - \kappa)$$

Now, we always consider that sparsity $s_0 = n^{\beta}$ with $0 < \beta < 1$. The dimensionality $p$ can either be a power of $n$, i.e., $p = n^k$ for some $k > 1$, or be a subexponential case whereby $n = \exp(n^{\gamma})$ where $0 < \gamma < 1$. In the second case, it is necessary to require that $\gamma + \beta < 1$ to ensure that $\frac{s_0}{n} \log\left(\frac{p}{s_0}\right) \to 0$ as $n \to \infty$. In both cases, the following result is true.

**Theorem 2.2.** *Given $n$ i.i.d. samples from (1.2), the minimax risk for estimating $\beta_0$ can be bounded from below as*

$$R_{\beta_0}(s_0, \alpha, \gamma) \gtrsim \frac{s_0}{n} \log\left(\frac{p}{s_0}\right); \tag{2.2}$$

*the minimax risk for estimating $f_0$ can be bounded from below as*

$$R_{f_0}(s_0, \alpha, \gamma) \gtrsim \max\left(n^{-2\alpha/2\alpha+1}, \frac{s_0}{n} \log\left(\frac{p}{s_0}\right)\right); \tag{2.3}$$

*moreover, the analogous result is also true for $R_{g_0}(s_0, \alpha, \gamma)$, namely*

$$R_{g_0}(s_0, \alpha, \gamma) \gtrsim \max\left(n^{-2\gamma/2\gamma+1}, \frac{s_0}{n} \log\left(\frac{p}{s_0}\right)\right), \tag{2.4}$$

*respectively.*

As is common in the literature, minimax lower bounds are obtained under the Gaussianity assumption on both errors $\varepsilon$ and the design matrix $X$. Such assumptions are meant to use known results on functional distances between normal density functions; e.g., Verzelen (2012).

As discussed previously, these lower bound results indicate (i) the best possible estimation of $\beta_0$ is not affected by the existence of nonparametric components, and coincides with the sparse estimation rate in high dimensional linear models; (ii) (the first order) minimax risk for estimating one nonparametric component does not depend on the smoothness of another component, but on the dimensionality and sparsity of the Euclidean parameter; see Figure 1. A similar lower bound has been discovered in nonparametric additive models (Raskutti et al., 2012) for the entire conditional mean function $\sum_{j \in S} h_{j0}(W_j)$, but with rather different interpretation: the term $s_0 \log(p/s_0)/n$ reflects the difficulty of selecting the sample size needed to perform the subset selection. Rather, this term here reflects the difficulty of selecting the $p$-dimensional vector $\beta_0$ with $s_0$ sparsity.

## 2.2. Nearly Optimal Estimators

In this section, we demonstrate that the penalized least square estimate for $(\beta_0, f_0, g_0)$ can almost achieve the lower bounds established in Theorem 2.2. To show such a result, we develop a series of oracle inequalities that give more and more refined estimation errors for each model component in the order of $g, f$ and $\beta$ (under the assumption that $f$ is smoother than $g$), and then derive the risk upper bounds by strengthening these oracle inequalities to their moment versions. Similar proof strategy was adopted in Müller and van de Geer (2015) and van de Geer and Muro (2015) to show oracle rates for parameters under partial linear models and nonparametric additive models, respectively. In comparison with Müller and van de Geer (2015), our nonparametric part possess an additive structure, and the linear covariates are relaxed from being bounded to sub-Gaussian.

Let $(\widehat{\beta}, \widehat{f}, \widehat{g})$ be an estimator of $(\beta_0, f_0, g_0)$ as follows:

$$(\widehat{\beta}, \widehat{f}, \widehat{g}) = \mathop{\arg\min}_{\beta \in \mathbb{R}^p, f \in W^{\alpha,2}(L_1), g \in W^{\gamma,2}(L_2)} \{\|Y - X^T\beta - f - g\|_n^2 + \lambda\|\beta\|_1 + \rho^2 J_\alpha^2(f) + \mu^2 J_\gamma^2(g)\}, \tag{2.5}$$

Without loss of generality, we assume that $\alpha \geq \gamma$.

**Assumption A.1.** *The covariates $X$ is a sub-Gaussian vector such that for any vector $v \in \mathbb{R}^p$, $v^T X$ is sub-Gaussian. And it satisfies for some constant $K_X \geq 1$,*

$$\sup_{v \in \mathbb{R}^p : \|v\|=1} \|v^T X\|_\Psi \leq K_X,$$

*where $\|\cdot\|_\Psi := \inf\{L > 0 : \mathbb{E}\Psi(\xi/L) < 1\}$ with $\Psi(t) = \exp(t^2) - 1$ is the Orlicz norm.*

**Assumption A.2.** *The error term $\varepsilon$ is independent of $(X, Z, U)$, and satisfies for some constant $K_\varepsilon \geq 1$,*

$$\|\varepsilon\|_\Psi \leq K_\varepsilon.$$

Let $\mathcal{H} = W^{\alpha,2}(L_1) \oplus W^{\gamma,2}(L_2)$ be a space of additive functions. For each $1 \leq j \leq p$, define the projection of $X_j$ onto $\mathcal{H}$ as $\Pi(X_j|\mathcal{H})(Z,U) = \arg\min_{h^* \in \mathcal{H}} \|X_j - h^*\|^2$. For simplicity, we write $((\Pi(X_1|\mathcal{H})(Z,U), \cdots, (\Pi(X_p|\mathcal{H})(Z,U))^T$ as $\pi_{X|Z,U}$. Note that $\pi_{X|Z,U} \in \mathbb{R}^p$ can be written as a sum of $f_X(Z) + g_X(U)$ where $f_{X_j} \in W^{\alpha,2}(L_1)$ and $g_{X_j} \in W^{\gamma,2}(L_2)$ for $1 \leq j \leq p$. Further, we have the following useful decomposition:

$$\|X^T\beta + f + g\|^2 = \|\widetilde{X}^T\beta\|^2 + \|\pi_{X|Z,U}^T\beta + f + g\|^2, \tag{2.6}$$

where $\widetilde{X} = X - \pi_{X|Z,U}$ is a random vector in $\mathbb{R}^p$.

**Assumption A.3.** *The smallest eigenvalue $\Lambda_{\min}^2$ of $\mathbb{E}\widetilde{X}\widetilde{X}^T$ is positive, and the largest eigenvalue $\Lambda_{\max}^2$ of $\mathbb{E}\{\pi_{X|Z,U}\pi_{X|Z,U}^T\}$ is finite.*

Assumption A.3 is common in semiparametric literature, e.g., Yu et al. (2011); Müller and van de Geer (2015). It guarantees that

$$\|\beta_{S_0}\|_1^2 \leq (\beta^T \mathbb{E}\widetilde{X}^T\widetilde{X}\beta)s_0/\Lambda_{\min}^2.$$

Our next assumption implies separate rates for $f$ and $g$ from that for $f+g$. This is due to $\|f+g\|^2 \geq (1-\gamma_0)(\|f\|+\|g\|)^2$, given $\mathbb{E}f_0(Z) = 0$, see Lemma 5.1 of van de Geer and Muro (2015) . Here, $\gamma_0$ is related to the minimal angle between two Hilbert spaces $W^{\alpha,2}(L_1)$ and $W^{\gamma,2}(L_2)$, see A.4 of Bickel et al. (1993), and formally defined as follows

$$\gamma_0^2 = \int (r-1)^2 p_Z p_U d\nu,$$

where $p = d\mathbb{P}_{ZU}/d\nu$ is the density of $\mathbb{P}_{ZU}$ w.r.t. $\nu = \nu_Z \times \nu_U$ with marginal densities $p_Z$ and $p_U$, and $r(z,u) = p(z,u)/(p_Z(z)p_U(u))$.

**Assumption A.4.** *It holds that $\gamma_0 < 1$.*

We assume the projection $f_P(U) = \mathbb{E}(f(Z)|U)$ to be smooth.

**Assumption A.5.** *For some constant $\Gamma > 0$, it holds that, for any function $f \in W^{\alpha,2}(L_1)$,*

$$J_\gamma(f_P) \leq \Gamma \|f\|.$$

**Remark.** Assumptions A.1 requires $X$ being sub-Gaussian, which relaxes the assumption that the entries of $X$ are uniformly bounded (say, by $M > 0$) in Müller and van de Geer (2015). In Müller and van de Geer (2015), the authors derive an upper bound on $\sup_{\beta,f} |(\mathbb{P}_n - \mathbb{P})X^T\beta f|$ in terms of $M$, and no other distributional information of $X$ is needed in this upper bound (as far as we are aware of). Our approach is more refined in that we bound $\sup_{\beta,f} |(\mathbb{P}_n - \mathbb{P})X^T\beta f|$ in terms of the Orlicz norm of the entries of $X$ (see Assumption A.1), which in general can be much smaller than $M$ even when $X$ is bounded. The relaxation of $X$ from being bounded to sub-Gaussian also leads to a more refined upper bound for $\sup_\beta |\mathbb{P}_n\varepsilon X^T\beta|$. Such a relaxation on $X$ is also needed in proving the minimax lower bound. Assumption A.2 relaxes the errors from being standard normal to sub-Gaussian, compared to Müller and van de Geer (2015). Condition 2.2 of van de Geer and Muro (2015) is not assumed as it holds up to a constant under our setting, see Lemma A.3 in the Appendix.

Before presenting our second main theorem, we need a set of oracle inequalities that hold in probability. Define the norm

$$\begin{aligned}
\tau(\beta, f, g; R) &= \lambda\frac{\|\beta\|_1}{\delta_0 R} + \|X^T\beta + f + g\| + \rho J_\alpha(f) + \mu J_\gamma(g), \\
\tau_I(\beta, f; R_I) &= \lambda\frac{\|\beta\|_1}{\delta_I R_I} + \|\widetilde{X}^T\beta\| + \|f_X^T\beta + f\| + \rho J_\alpha(f),
\end{aligned}$$

for some constant $\delta_0 > 0$.

**Lemma 2.3.** *Suppose Assumptions A.1-A.5 hold. Let*

$$\lambda \gtrsim \sqrt{\frac{\log p}{n}}, \quad \rho^2 \asymp n^{-\frac{2\alpha}{2\alpha+1}}, \quad \mu^2 \asymp n^{-\frac{2\gamma}{2\gamma+1}} \quad \textit{and} \quad \rho^2 \leq \mu^2.$$

*If there exist $R$ and $R_I$ satisfying*

$$\frac{4\lambda^2 s_0}{\Lambda^2_{\min}} \leq R_I^2 \leq R^2 \leq \lambda \leq 1, \quad \rho^2 \leq \frac{\delta_I^2 R_I^2}{2(1+\Gamma+L_1+J_\alpha(f_0))^2}, \quad \mu^2 \leq \frac{\delta_0^2 R^2}{2(1+\Gamma+L_2+J_\gamma(g_0))^2},$$

*then it holds that*

$$\mathbb{P}\left(\tau(\widehat{\beta}-\beta_0, \widehat{f}-f_0, \widehat{g}-g_0; R) \leq R, \tau_I(\widehat{\beta}-\beta_0, \widehat{f}-f_0; R_I) \leq R_I\right) \geq 1 - C\exp(-n\rho^2/c)$$

*for some constants $C, c > 0$.*

In particular, we can take $R^2 \asymp \mu^2 + \lambda^2 s_0$ and $R_I^2 \asymp \rho^2 + \lambda^2 s_0$. Then the first oracle inequality gives an upper bound for the overall estimating rate of $(\widehat{\beta}, \widehat{f}, \widehat{g})$:

$$O_P\left(\max\left(n^{-2\gamma/(2\gamma+1)}, s_0 \log p/n\right)\right),$$

which implies the desirable estimation rate for $\widehat{g}$. And the second one provides a tighter bound for the estimating rate of $(\widehat{\beta}, \widehat{f})$:

$$O_P\left(\max\left(n^{-2\alpha/(2\alpha+1)}, s_0 \log p/n\right)\right),$$

which in turn implies the rate for $\widehat{f}$.

We need a separate lemma to improve the rate of $\|\widehat{\beta}-\beta_0\|$ to (nearly) minimax optimal level $s_0 \log(p)/n$. This new Lemma 2.4 requires us to project $X$ onto the additive space $\mathcal{H}$.

**Lemma 2.4.** *Assume conditions of Lemma 2.3 hold. Then there exists constants $C', c' > 0$ such that with probability at least $1 - 7/(2p) - C'\exp(-c'n\rho^2)$,*

$$\|\widetilde{X}^T(\widehat{\beta} - \beta_0)\|_n^2 + (\lambda/2)\|\widehat{\beta} - \beta_0\|_1 \leq \frac{4s_0\lambda^2}{\Lambda_{\min}^2}.$$

Lemma 2.4 has two important implications: (i) prediction error: $\|\widetilde{X}^T(\widehat{\beta} - \beta_0)\|_n^2 \leq 4s_0\lambda^2/\Lambda_{\min}^2$; (ii) $\ell_1$ error: $\|\widehat{\beta} - \beta_0\|_1 \leq 8s_0\lambda/\Lambda_{\min}^2$. We note that these two rates are in the same order as those standard lasso rates (as if $f_0$ and $g_0$ were known); see Bühlmann and van de Geer (2011). However, the probability that these rates hold is comparatively smaller as reflected by an additional term $\exp(-c'n\rho^2)$. This is the price to pay for estimating two unknown nonparametric functions in the model.

We are now ready to prove that $(\widehat{\beta}, \widehat{f}, \widehat{g})$ nearly achieve the minimax lower bounds established in Theorem 2.2.

**Theorem 2.5.** *Assume conditions of Lemma 2.3 hold. Then*

$$\mathbb{E}\|\widehat{\beta} - \beta_0\|^2 \lesssim \frac{s_0 \log p}{n},$$

$$\mathbb{E}\int_0^1 |\widehat{f}(z) - f_0(z)|^2\, dz \lesssim \max\left(n^{-2\alpha/(2\alpha+1)}, \frac{s_0 \log p}{n}\right),$$

*and*

$$\mathbb{E}\int_0^1 |\widehat{g}_0(u) - g_0(u)|^2\, du \lesssim \max\left(n^{-2\gamma/(2\gamma+1)}, \frac{s_0 \log p}{n}\right).$$

# APPENDIX

In this section, we present all the technical details. Proofs for main lemmas, theorems and corollaries in Sections 2.1, 2.2 are presented in Sections A.1, A.2, respectively. Results from empirical process theory are presented in Section A.3.

## A.1. Proofs for Section 2.1

### *A.1.1. Proof of Theorem 2.2*

***Proof.*** It is easy to see that the minimax lower bound for estimating $\beta_0$ trivially follows from that for high dimensional linear models derived in Verzelen (2012), i.e., $R_{\beta_0}(s_0, \Sigma, \alpha, \gamma) \geq \inf_{\widehat{\beta}} \sup_{\beta_0 \in B[s_0,p]} \mathbb{E}[\|\beta_0 - \widehat{\beta}\|^2]$ (fixing $f$ and $g$ at their true values).

In what follows, we concentrate on the lower bound of the minimax risk for estimators of the nonparametric component. Without loss of generality, we choose $f_0$ for our discussion. To make this proof easier, we start from partial linear models

$$Y = X^T\beta_0 + f_0(Z) + \varepsilon, \tag{A.1}$$

where $\beta_0 \in B[s_0, p]$ and $f_0 \in W^{\alpha,2}(L_1)$, and will show the minimax risk for $f_0$ is bounded from below by

$$\max\left(n^{-2\alpha/(2\alpha+1)}, \frac{s_0}{n}\log\left(\frac{p}{s_0}\right)\right) \tag{A.2}$$

up to a universal constant, based on iid observations $\{Y_i, X_i, Z_i\}_{i=1}^n$.

In the model (A.1), define the minimax estimation risk for $f_0$ as

$$R_{f_0}(s_0, \alpha) = \inf_{\Sigma \in \mathcal{S}_p} \inf_{\widehat{f}} \sup_{\beta_0 \in B[s_0,p]} \sup_{f_0 \in W^{\alpha,2}(L_1)} \mathbb{E}\int_0^1 |\widehat{f}(z) - f_0(z)|^2\, dz.$$

The first part of lower bound, i.e., $n^{-2\alpha/(2\alpha+1)}$, trivially follows from the following inequality (assuming $\beta$ taking its true value)

$$R_{f_0}(s_0,\alpha) \geq \inf_{\Sigma\in\mathcal{S}_p} \inf_{\widehat{f}} \sup_{f_0\in W^{\alpha,2}(L_1)} \mathbb{E}\int_0^1 |\widehat{f}(z) - f_0(z)|^2\, dz$$

and the classical nonparametric minimax rate.

Our method of obtaining minimax lower bounds on rates of convergence for estimators of the nonparametric component is somewhat different from typical ones. It is based on the the Corollary 2.19 from Massart (2007), that represents a version of the classical Fano's lemma. Specifically, Massart (2007) specifies that values of the unknown parameter are viewed as points in a pseudometric space $(S,d)$ where $d$ is the corresponding pseudometric. Note that both the coefficient $\beta$ and the function $f$ are not known yet it is only the function $f$ that is currently a quantity of interest. To help us handle this situation we define a "real" pseudometric (the one that is not a metric) between the two points of interest and use it to establish the lower bound. In this situation, defining the pseudometric between pairs $(\beta_1, f_1)$ and $(\beta_2, f_2)$ as the $L_2$-distance between $f_1$ and $f_2$ works out well.

Now we need to establish the second part of the minimax lower bound, i.e., $s_0\log(p/s_0)/n$. Using the approach just described, we note that for two vectors $s_1 = (\beta_1, f_1)^{'}$ and $s_2 = (\beta_2, f_2)^{'}$ we can define a pseudometric $d(s_1, s_2) := d_1(f_1, f_2)$. It is easy to verify that all of the metric properties are satisfied for $d(s_1, s_2)$ except that, of course, it is possible to have $d(s_1, s_2) = 0$ while $s_1 \neq s_2$; this, clearly, qualifies $d$ as a pseudometric. Choosing $r = 2$ reduces the search for the lower bound of $\sup_{s\in\mathcal{C}} \mathbb{E}\,[d^2(s,\widehat{s})]$ to that of $\sup_{s\in\mathcal{C}} \mathbb{E}\,[d_1^2(f,\widehat{f})]$.

Our first step is, thus, to find

$$\delta := \min_{\substack{s_1,s_2\in\mathcal{C}\\ s_1\neq s_2}} d_1(f_1, f_2)$$

for an appropriate $\mathcal{C}$. To define the set $\mathcal{C}$, we start with selecting a set of test functions $f_{jn}$ and $f_{kn}$ (note that they depend on $n$). To do so, first define a kernel function $K_0(u) = \exp\left(-\frac{1}{1-u^2}\right) I(|u|\leq 1)$ and take $K(u) = aK_0(2u)$ for a sufficiently small constant $a > 0$. For an integer $m \geq 1$, and $k = 1,\ldots,m$, select a set of points $z_k = \frac{k-1/2}{m}$ that belong in $[0,1]$. For convenience, we will also use the following notation: $\Delta_0 = [0, 1/m]$ and $\Delta_k = [(k-1)/m, k/m]$, where $k = 2,\ldots,m$. The choice of $m$ will depend on $n$, $s_0$, $p$ and $\alpha$. For brevity, we introduce the notation $\delta_n := \frac{s_0}{n}\log\frac{p}{s_0}$ where $\delta_n \to 0$ as $n\to\infty$. Now we can define the cardinality of the partition of $[0,1]$ as $m = n^{2\alpha+1/(4\alpha+1)}\delta_n^{2\alpha/(4\alpha+1)}$; the corresponding optimal bandwidth is defined as $h_n = \left(\frac{\delta_n}{m}\right)^{1/2\alpha+1}$. Note that the choices of $h_n$ and $m$ that we made are sensible since one can easily show that $m\to\infty$ as $n\to\infty$ while the bandwidth $h_n$ goes to zero as $n\to\infty$ as well. This is true both in the polynomial setting where $p = n^k$ and in subexponential setting where $p = \exp(n^\gamma)$. With all of the elements in place, we can now define a function $\phi_k(z) = L_1 h_n^\alpha K\left(\frac{z-z_k}{h_n}\right)$.

Secondly, consider a set of binary sequences $\Omega = \{\omega = (\omega_1,\ldots,\omega_m), \omega_i\in\{0,1\}\}$, and define a set of functions $\mathcal{F} = \{f_\omega(z) = \sum_{k=1}^m \omega_k\phi_k(z), \omega\in\Omega\}$. First of all, we note that any function $f_\omega\in\mathcal{F}$ belongs, by construction, to $W^{\alpha,2}(L_1)$. Second, we need to select test functions from the set $\mathcal{F}$; however, in practice, in order to ensure that any two functions thus selected are separated by at least the required amount, a certain *subset* of the set $\mathcal{F}$ has to be used. Due to the Varshamov-Gilbert lemma (see e.g. Gilbert (1952)) we can, indeed, find a subset of $\mathcal{F}$, i.e., $\{f_{jn}, j = 0,\ldots,M\}$, such that any two distinct functions in it are sufficiently well separated and whose cardinality $M$ is sufficiently large. More specifically, for any $0\leq j < k \leq M$, we have $f_{jn}$ and $f_{kn}$ such that the squared $L^2$ distance between the two functions is $d_1^2(f_{jn}, f_{kn}) \asymp \delta_n$ as long as $\log M \geq \frac{\log 2}{8}m$. This is done using a standard construction one can find in, for example, Tsybakov (2008). To show that this is true, denote the binary sequences corresponding to $f_{jn}$ and $f_{kn}$ as $\omega^{(j)}$ and $\omega^{(1/2\alpha)}$, respectively, while the Hamming distance between them is denoted $\rho(\omega^{(j)},\omega^{(1/2\alpha)})$. By the Varshamov-Gilbert lemma, we have

$$||f_{jn} - f_{kn}||_2^2 = L_1^2 h_n^{2\alpha+1}||K||_2^2 \rho(\omega^{(j)},\omega^{(1/2\alpha)}) \geq L_1^2 h_n^{2\alpha+1}||K||_2^2\frac{m}{16} \geq C\delta_n.$$

Define a finite set $\mathcal{C}$ that consists of $s = (\beta, f_{jn})'$ where $\beta$ is arbitrary while $0 \le j \le M$; in other words, the set $\mathcal{C}$ consists of all vectors with an arbitrary parametric first coordinate and one of the test functions we constructed as a second coordinate. The vectors thus constructed are all distinct so, by Fano's lemma, and using the pseudometric $d$, we immediately obtain that $\sup_{s\in\mathcal{C}} \mathbb{E}_s[d^2(s,\widehat{s})] = \sup_{s\in\mathcal{C}} \mathbb{E}_s[d_1^2(f,\widehat{f})] \ge C\delta_n$ as needed. Now it only remains to verify that $\max_{s,t\in\mathcal{C}} K(P_s,P_t) \le \kappa \log|C|$. To do so, we first note that the cardinality of $\mathcal{C}$, as defined, is $M$; using calculations very similar to those in Tsybakov (2008) (p. 115-116) we find that, for fixed data points $Z_1,\dots,Z_n$, we have $\max_{s,t\in\mathcal{C}} K(P_s,P_t) \le Cnh_n^{2\alpha} \le m$ for the $h_n$ we defined above. By Varshamov-Gilbert inequality, it follows that $m \le \frac{8}{\log 2}\log M$ and so the condition on the Kullback-Leibler distance is satisfied as well. Thus, the statement has been proved for the partial linear model.

Note that an interesting feature of our proof is the subtle way in which the bandwidth of the test functions $\phi_k(z)$ and the cardinality of the set of these functions $m$ depend on each other. The bandwidth $h_n = (\delta_n/m)^{1/2\alpha+1}$ and $m = n^{2\alpha+1/(4\alpha+1)}\delta_n^{2\alpha/(4\alpha+1)}$ where $\delta_n = (s_0/n)\log(p/s_0)$ guarantees the existence of the non-trivial lower bound due to Varshamov - Gilbert Lemma. Moreover, it works for a wide range of dimensionalities $p$ that includes both the polynomial setting $p = n^k$ and subexponential $p = \exp(n^\gamma)$ with $k,\gamma > 0$. The precise selection of the relationship between $m$ and $h_n$ is what enables us to obtain the correct lower bound of the risk.

To carry these results over to the partial linear additive model (1.2), we need to consider a nonparametric model without the linear component

$$Y_i = A + f_0(Z_i) + g_0(U_i) + \varepsilon_i, \tag{A.3}$$

where $A$ is a constant, $\mathbb{E} f_0(Z) = \mathbb{E} g_0(U) = 0$ for identifiability purposes, $(Z,U) \in [0,1]\times[0,1]$, and $f_0 \in W^{\alpha,2}(L_1)$ and $g_0 \in W^{\gamma,2}(L_2)$. For the model (A.3), it is known (see e.g. Horowitz et al. (2006)) that the minimax risk of estimating $f_0$ is $n^{-2\alpha/(2\alpha+1)}$, which does not depend on $\gamma$. Specifically, this means that

$$\inf_{\widehat{f}} \sup_{g_0\in W^{\gamma,2}(L_2)} \sup_{f_0\in W^{\alpha,2}(L_1)} \mathbb{E}\int_0^1 |\widehat{f}(z) - f_0(z)|^2\,dz \ge Cn^{-2\alpha/(2\alpha+1)}$$

for a generic constant $C$ that does not depend on $n$. By the definition of $R_{f_0}(s_0,\alpha,\gamma)$, this immediately suggests that one lower bound of $R_{f_0}(s_0,\alpha,\gamma)$ is $n^{-2\alpha/(2\alpha+1)}$. On the other hand, it is also clear that (assuming $g$ at its true value)

$$R_{f_0}(s_0,\alpha,\gamma) \ge \inf_{\widehat{f}} \sup_{\beta_0\in B[s_0,p]} \sup_{f_0\in W^{\alpha,2}(L_1)} \mathbb{E}\int_0^1 |\widehat{f}(z) - f_0(z)|^2\,dz.$$

Then, by the lower bound result for partial linear models, we know that (A.2) is another lower bound for estimating $f$ in partial linear additive models. This concludes our proof. □

## A.2. Proof for Section 2.2

### *A.2.1. Proof of Lemma 2.3*

Before proving the Lemma, we first present the following necessary notations.

For any normed linear space $\mathcal{F}$, let $d$ be a metric on the space $\mathcal{F}$. For any $t > 0$, define $N(t,\mathcal{F},d)$ as covering number of $\mathcal{F}$ and $H(t,\mathcal{F},d) = \log N(t,\mathcal{F},d)$ as entropy number of $\mathcal{F}$. Let $\mathcal{A}_n$ be the set of all configurations $A_n$ of $n$ points within the support of the joint density $P_{XZU}$. For $A_n \in \mathcal{A}_n$, $\|f\|_{A_n,\infty} := \max_{Z\in A_n}|f(Z)|$. Let $\mathcal{H}_\infty(t,\mathcal{F}) = \sup_{A_n\in\mathcal{A}_n} H(t,\mathcal{F},\|\cdot\|_{A_n,\infty})$, see van de Geer (2014). Further, we write

$$\mathcal{J}_\infty(u,\mathcal{F}) = C_0 \inf_{\delta>0}\left[u\int_{\delta/4}^1 \sqrt{\mathcal{H}_\infty(tu/2,\mathcal{F})}dt + \sqrt{n}\delta u\right].$$

For arbitrary constants $R_0 > 0$ and $M_0 > 0$, we denote $W^{\alpha,2}(R_0,M_0) = \{f \in W^{\alpha,2}(L_1) : \|f\| \le R_0, J_\alpha(f) \le M_0\}$ and $W^{\gamma,2}(R_0,M_0) = \{g \in W^{\gamma,2}(L_2) : \|g\| \le R_0, J_\gamma(g) \le M_0\}$. Therefore, it holds

that for $R_0 \leq M_0$ and some constants $A_I \geq 1$ and $A_J \geq 1$

$$\mathcal{J}_\infty(z, W^{\alpha,2}(R_0, M_0)) \leq A_I M_0^{1/2\alpha} z^{1-1/2\alpha}, \tag{A.4}$$

and

$$\mathcal{J}_\infty(z, W^{\gamma,2}(R_0, M_0)) \leq A_J M_0^{1/2\gamma} z^{1-1/2\gamma}.$$

For some $\delta_0 > 0$ small enough, define

$$\mathcal{M}(R) = \{(\beta, f, g) : \tau(\beta, f, g; R) \leq R, \beta \in \mathbb{R}^p, f \in W^{\alpha,2}(L_1), g \in W^{\gamma,2}(L_2)\},$$

$$\mathcal{T}_1(R) = \left\{ \sup_{\mathcal{M}(R)} \left| \|X^T\beta + f + g\|_n^2 - \|X^T\beta + f + g\ \|^2 \ \right| \leq \delta_0^2 R^2 \right\},$$

$$\mathcal{T}_2(R) = \left\{ \sup_{\mathcal{M}(R)} \left| \mathbb{P}_n\left(\varepsilon(X^T\beta + f + g)\right) \right| \leq \delta_0^2 R^2 \right\},$$

and

$$\mathcal{T}(R) = \mathcal{T}_1(R) \cap \mathcal{T}_2(R). \tag{A.5}$$

Let $f_{XP}(\cdot) = \mathbb{E}(f_X(Z)|U = \cdot)$ and $f_{XA} = f_X - f_{XP}$. For any $f \in W^{\alpha,2}(L_1)$, write $f_P(\cdot) = \mathbb{E}(f(Z)|U = \cdot)$ and $f_A = f - f_P$. For $\delta_I$ sufficiently small, define

$$\mathcal{M}_I(R_I) = \left\{(\beta, f) : \tau_I(\beta, f; R_I) \leq R_I, \beta \in \mathbb{R}^p, f \in W^{\alpha,2}(L_1)\right\},$$

$$\mathcal{T}_{I,1}(R_I) = \left\{ \sup_{(\beta,f)\in\mathcal{M}_I(R_I)} \left| \|\widetilde{X}^T\beta + f_{XA}^T\beta + f_A\|_n^2 - \|\widetilde{X}^T\beta + f_X^T\beta + f_A\|^2 \right| \leq \delta_I^2 R_I^2 \right\},$$

$$\mathcal{T}_{I,2}(R_I) = \left\{ \sup_{(\beta,f)\in\mathcal{M}_I(R_I)} \left| \mathbb{P}_n\left(\varepsilon(\widetilde{X}^T\beta + f_{XA}^T\beta + f_A)\right) \right| \leq \delta_I^2 R_I^2 \right\},$$

$$\mathcal{T}_{I,3}(R_I) = \left\{ \sup_{(\beta,f,g)\in\mathcal{M}(R),(\beta,f)\in\mathcal{M}_I(R_I)} \left| \mathbb{P}_n(\widetilde{X}^T\beta + f_{XA}^T\beta + f_A)(f_{XP}^T\beta + g_X^T\beta + f_P + g) \right| \leq \delta_I^2 R_I^2 \right\},$$

and

$$\mathcal{T}_I(R_I) = \mathcal{T}_{I,1}(R_I) \cap \mathcal{T}_{I,2}(R_I) \cap \mathcal{T}_{I,3}(R_I). \tag{A.6}$$

To prove Lemma 2.3, we first show in Lemma A.1 that $\tau(\widehat{\beta} - \beta_0, \widehat{f} - f_0, \widehat{g} - g_0; R) \leq R$ on $\mathcal{T}(R)$. The probability of $\mathcal{T}(R)$ is estimated in Lemma A.4. We next show $\tau_I(\widehat{\beta} - \beta_0, \widehat{f} - f_0; R_I) \leq R_I$ on the set $\mathcal{T}(R) \cap \mathcal{T}_I(R_I)$ in Lemma A.5, whereas the probability of $\mathcal{T}_I(R_I)$ is estimated in Lemma A.6. Lemmas A.2 and A.3 are technical Lemmas in order to show Lemmas A.4 and A.6.

**Lemma A.1.** *Under the conditions of Lemma 2.3, we have, on $\mathcal{T}(R)$,*

$$\tau(\widehat{\beta} - \beta_0, \widehat{f} - f_0, \widehat{g} - g_0; R) \leq R.$$

***Proof.*** Take $\delta_0 \leq 1/30$. Under the conditions of Lemma 2.3, we can find $\rho$ and $\mu$ such that

$$\rho^2 J_\alpha^2(f_0) + \mu^2 J_\gamma^2(g_0) \leq \delta_0^2 R^2, \tag{A.7}$$

and

$$4\lambda^2 s_0/\Lambda_{\min}^2 \leq R_I^2 \leq R^2. \tag{A.8}$$

Define

$$t = \frac{R}{R + \tau(\widehat{\beta} - \beta_0, \widehat{f} - f_0, \widehat{g} - g_0; R)}.$$

Let $\widetilde{\beta} = t\widehat{\beta} + (1-t)\beta_0$, $\widetilde{f} = t\widehat{f} + (1-t)f_0$, $\widetilde{g} = t\widehat{g} + (1-t)g_0$. Notice that $\tau(\widetilde{\beta} - \beta_0, \widetilde{f} - f_0, \widetilde{g} - g_0; R) = t\tau(\widehat{\beta} - \beta_0, \widehat{f} - f_0, \widehat{g} - g_0; R) \le R$, which implies $(\widetilde{\beta} - \beta_0, \widetilde{f} - f_0, \widetilde{g} - g_0) \in \mathcal{M}(R)$. In order to show $\tau(\widehat{\beta} - \beta_0, \widehat{f} - f_0, \widehat{g} - g_0; R) \le R$, it suffices to prove $\tau(\widetilde{\beta} - \beta_0, \widetilde{f} - f_0, \widetilde{g} - g_0; R) \le R/2$.

By the convexity, we have

$$\begin{aligned}&\|Y - X^T\widetilde{\beta} - \widetilde{f} - \widetilde{g}\|_n^2 + \lambda\|\widetilde{\beta}\|_1 + \rho^2 J_\alpha^2(\widetilde{f}) + \mu^2 J_\gamma^2(\widetilde{g})\\ \le&\|Y - X^T\beta_0 - f_0 - g_0\|_n^2 + \lambda\|\beta_0\|_1 + \rho^2 J_\alpha^2(f_0) + \mu^2 J_\gamma^2(g_0).\end{aligned}$$

Together with (A.7), it further implies

$$\begin{aligned}&\|X^T(\widetilde{\beta} - \beta_0) + (\widetilde{f} - f_0) + (\widetilde{g} - g_0)\|_n^2 + \lambda\|\widetilde{\beta}\|_1 + \rho^2 J_\alpha^2(\widetilde{f}) + \mu^2 J_\gamma^2(\widetilde{g})\\ \le&2\mathbb{P}_n\left(\varepsilon(X^T(\widehat{\beta} - \beta_0) + (\widehat{f} - f_0) + (\widehat{g} - g_0))\right) + \lambda\|\beta_0\|_1 + \delta_0^2 R^2. \qquad \text{(A.9)}\end{aligned}$$

Therefore, by the definition of $\mathcal{T}_1(R)$ and $\mathcal{T}_2(R)$,

$$\begin{aligned}&\|X^T(\widetilde{\beta} - \beta_0) + (\widetilde{f} - f_0) + (\widetilde{g} - g_0)\|^2 + \lambda\|\widetilde{\beta}_{S_0^C}\|_1 + \rho^2 J_\alpha^2(\widetilde{f}) + \mu^2 J_\gamma^2(\widetilde{g})\\ \le&\delta_0^2 R^2 + \delta_0^2 R^2 + 2\delta_0^2 R^2 + \lambda\|\beta_0\|_1 - \lambda\|\widetilde{\beta}_{S_0}\|_1\\ \le&4\delta_0^2 R^2 + \lambda\|\beta_{0S_0} - \widetilde{\beta}_{S_0}\|_1. \qquad \text{(A.10)}\end{aligned}$$

Note that

$$\begin{aligned}\lambda\|\beta_{0S_0} - \widetilde{\beta}_{S_0}\|_1 &\le \lambda\sqrt{s_0}\|\beta_{0S_0} - \widetilde{\beta}_{S_0}\|\\ &\le \lambda\sqrt{s_0}\|\widetilde{\beta} - \beta_0\|\\ &\le \lambda\sqrt{s_0}\|\widetilde{X}^T(\widetilde{\beta} - \beta_0)\|/\Lambda_{\min}\\ &\le \lambda^2 s_0/\Lambda_{\min}^2 + \|\widetilde{X}^T(\widetilde{\beta} - \beta_0)\|^2/4\\ &\le \delta_0^2 R^2/4 + \|\widetilde{X}^T(\widetilde{\beta} - \beta_0)\|^2/4, \qquad \text{(A.11)}\end{aligned}$$

where the third inequality holds by Assumption A.3, the fourth inequality follows from $uv \le u^2 + v^2/4$, and the last one is due to (A.8). Thus, substituting (A.11) into (A.10), we obtain

(a) $(3/4)\|X^T(\widetilde{\beta} - \beta_0) + (\widetilde{f} - f_0) + (\widetilde{g} - g_0)\|^2 \le (17/4)\delta_0^2 R^2$, by orthogonal decomposition (2.6);
(b) $\rho^2 J_\alpha^2(\widetilde{f}) \le (17/4)\delta_0^2 R^2$;
(c) $\mu^2 J_\gamma^2(\widetilde{g}) \le (17/4)\delta_0^2 R^2$.

Now it follows from (a) that $\|X^T(\widetilde{\beta} - \beta_0) + (\widetilde{f} - f_0) + (\widetilde{g} - g_0)\| \le (\sqrt{17}/\sqrt{3})\delta_0 R$. In addition, (b), (c) and (A.7) imply

$$\rho I(\widetilde{f} - f_0) \le \rho J_\alpha(\widetilde{f}) + \rho J_\alpha(f_0) \le \frac{\sqrt{17}}{2}\delta_0 R + 2\delta_0 R \le \sqrt{17}\delta_0 R$$

and

$$\left(\frac{\mu}{R}\right)^{\frac{2-q}{q}} J_\gamma^2(\widetilde{g} - g_0) \le \frac{\sqrt{17}}{2}\delta_0 R + 2\delta_0 R \le \sqrt{17}\delta_0 R.$$

Adding $\lambda\|\beta_{0S_0} - \widetilde{\beta}_{S_0}\|_1$ on both sides of (A.10), we get $\|X^T(\widetilde{\beta} - \beta_0) + (\widetilde{f} - f_0) + (\widetilde{g} - g_0)\|^2 + \lambda\|\widetilde{\beta} - \beta_0\|_1 + \rho^2 J_\alpha^2(\widetilde{f}) + \mu^2 J_\gamma^2(\widetilde{g}) \le 4\delta_0^2 R^2 + 2\lambda\|\beta_{0S_0} - \widetilde{\beta}_{S_0}\|_1 \le 4\delta_0^2 R^2 + \|\widetilde{X}^T(\widetilde{\beta} - \beta_0)\|^2 + \frac{1}{4}\delta_0^2 R^2$, which further implies

$$\lambda\|\widetilde{\beta} - \beta_0\|_1 \le \frac{17}{4}\delta_0^2 R^2.$$

Invoking the definition of $\tau(\cdot)$, we finally get

$$\tau(\widetilde{\beta} - \beta_0, \widetilde{f} - f_0, \widetilde{g} - g_0; R) \le \left(\sqrt{17}/\sqrt{3} + 2\sqrt{17} + 17/4\right)\delta_0 R \le 15\delta_0 R \le \frac{1}{2}R,$$

by letting $\delta_0 \le 1/30$. □

**Lemma A.2.** *Given $\xi_i, 1 \le i \le n$ are i.i.d. random sub-Gaussian vectors with covariance matrix $\Xi$*

$$\mathbb{P}\left(\sup_{\beta \in B(s_0,p)} \left(\frac{\beta^T \left(\frac{1}{n}\sum_{i=1}^n \xi_i \xi_i^T\right)\beta}{\beta^T \Xi \beta} - 1\right) > C\left(\sqrt{\frac{t + s_0 \log p}{n}} + \frac{t + s_0 \log p}{n}\right)\right) \le \exp(-t).$$

*where $C > 0$ is a constant not depending on $n$.*

This Lemma follows the same reasoning as Lemma 1 of Nickl and van de Geer (2013).

**Lemma A.3** (Gagliardo–Nirenberg–Sobolev inequality)**.** *For $\alpha \ge 1$, there exists a constant $C_\alpha \ge 1$ such that*

$$\sup_{f \in W^{\alpha,2}(R_0,M_0)} \|f\|_\infty \le C_\alpha (R_0 + M_0).$$

*Remark.* The inequality above is standard for Sobolev spaces consisting of functions that vanish at the endpoints 0 and 1, or of 1-periodic functions. In our paper, the Sobolev space $W^{m,2}(R_0, M_0)$ consists of functions that are (one-sided) differentiable at the endpoints (such as splines), therefore differs slightly from the commonly used definition. For such Sobolev space, we had not been able to locate a proof in the literature, therefore include here a proof for readers' convenience.

***Proof of Lemma A.3.*** Fix $x \in [0, 1]$. Using Taylor's theorem, we can write for any $t \in [0, 1]$

$$f(x) = f(t) + Q_{\alpha-1}(x, t) + R_\alpha(x, t)$$

where

$$Q_{\alpha-1}(x, t) = f'(t)(x - t) + \cdots + \frac{f^{(\alpha-1)}(t)}{(\alpha - 1)!}(x - t)^{\alpha-1}$$

and

$$R_\alpha(x, t) = \int_t^x \frac{f^{(\alpha)}(s)}{(\alpha - 1)!}(x - s)^{\alpha-1} ds$$

is the Lagrange remainder. Averaging over $t$, we obtain

$$f(x) = \int_0^1 f(t)dt + \int_0^1 Q_{\alpha-1}(x, t)dt + \int_0^1 R_\alpha(x, t)dt.$$

It is easy to see that

$$\left|\int_0^1 f(t)dt\right| \le \|f\| \le R_0, \quad \left|\int_0^1 R_m(x, t)dt\right| \le \|f^{(m)}\| \le M_0.$$

For the middle term, we have

$$\left|\int_0^1 Q_{m-1}(x, t)dt\right| \le \|f'\|_\infty + \cdots + \frac{\|f^{(m-1)}\|_\infty}{(m - 1)!}.$$

So it suffices to show

$$\sup_{f \in W^{m,2}(R_0,M_0)} \|f^{(k)}\|_\infty \le C_m(R_0 + M_0), \; 1 \le k \le m - 1. \tag{A.12}$$

The proof of (A.12) is by induction on $m - k = 1, \cdots, m - 1$. The base case is $k = m - 1$. By averaging over suitable $k$-th order finite differences of $f$, we have

$$\min_{x \in [0,1]} |f^{(k)}(x)| \le C_m \|f\| \le C_m R_0. \tag{A.13}$$

Suppose the minimum is attained at $x_0$, i.e.

$$|f^{(k)}(x_0)| = \min_{x\in[0,1]} |f^{(k)}(x)|.$$

Then, for any $x \in [0,1]$, we can write

$$f^{(k)}(x) = f^{(k)}(x_0) + \int_{x_0}^{x} f^{(k+1)}(t)dt.$$

Thus

$$|f^{(k)}(x)| \le |f^{(k)}(x_0)| + \|f^{(k+1)}\|.$$

Combining (A.13) and the bound on $\|f^{(k+1)}\|$, this implies

$$\|f^{(k)}\|_\infty \le C_m(R_0 + M_0).$$

By induction, the same bound holds for $k < m-1$, with the same argument. This establishes (A.12), and the proof of Lemma A.3 is complete. □

**Lemma A.4.** *Under the conditions of Lemma 2.3, we have for some constants* $\widetilde{C} > 0, \widetilde{c} > 0$, $\mathbb{P}(\mathcal{T}(R)) \ge 1 - \widetilde{C}\exp(-\widetilde{c}n\rho^2)$.

***Proof.*** We first introduce some necessary notations and preliminary results. Note that $\tau(\beta, f, g; R) \le R$ implies that $\|X^T\beta + f + g\|^2 \le R^2$ and $\|\beta\|_1 \le \delta_0 R^2/\lambda$, $J_\alpha(f) \le R/\rho$, $J_\gamma(g) \le R/\mu$. By the orthogonal decomposition (2.6), we have $\|\widetilde{X}^T\beta\| \le R$ and $\|f_X^T\beta + f + g_X^T\beta + g\| \le R$. Then Assumption A.3 implies

$$\|X^T\beta\| \le \|\widetilde{X}^T\beta\| + \|\pi^T_{X|Z,U}\beta\| \le R + (\Lambda_{\max}/\Lambda_{\min})\|\widetilde{X}^T\beta\| \le (1 + \Lambda_{\max}/\Lambda_{\min})R.$$

Similar arguments and Assumption A.4 imply that both $\|f\|$ and $\|g\|$ are bounded by $(1+\Lambda_{\max}/\Lambda_{\min})R/\sqrt{(1-\gamma_0)}$, i.e. $R/\sqrt{1-\gamma_1}$, for simplicity, we write it as $R_1$. Define $M_2 = R/\rho$, $M_3 = R/\mu$. In particular, we can choose $\rho^2 \le 1-\gamma_1$ and $\mu^2 \le 1-\gamma_1$, where $1-\gamma_1 = (1-\gamma_0)/(1+\Lambda_{\max}/\Lambda_{\min})^2$. Then it follows from Lemma A.3 that $\sup_{f\in W^{\alpha,2}(R_1,M_2)} \|f\|_\infty \le C_\alpha M_2$ and $\sup_{g\in W^{\gamma,2}(R_1,M_3)} \|g\|_\infty \le C_\gamma M_3$, with the fact $R_1 \le M_2$ and $R_1 \le M_3$. Further, we find a constant $L > 1$ such that the following hold:

$$\sqrt{n}\rho^{\frac{2\alpha+1}{2\alpha}} \ge LA_I,\ \sqrt{n}\mu^{\frac{2\gamma+1}{2\gamma}} \ge LA_J,\ R \ge \mu \ge \rho,\ \rho^{1/2\alpha} \le 1/L, \sqrt{s_0 \log p/n} \le 1/L. \tag{A.14}$$

Now, we are ready to apply empirical process theory stated in Section A.3 to show that with probability at least $1 - 6\exp(-n\rho^2/L)$, the event $\mathcal{T}_1(R)$ holds. Without loss of generality, we take $C_1 = 1$ in Theorem A.11; otherwise we can replace in $L = LC_1$ in the proof. Note that for any $(\beta, f, g)$, it holds

$$\begin{aligned}
&\big|\|X^T\beta + f + g\|_n^2 - \|X^T\beta + f + g\ \|^2\ \big| \\
\le&\big|\|X^T\beta\|_n^2 - \|X^T\beta\|^2\big| + \big|\|f\|_n^2 - \|f\|^2\big| + \big|\|g\|_n^2 - \|g\|^2\big| + \big|2(\mathbb{P}_n - \mathbb{P})X^T\beta f\big| \\
&+ \big|2(\mathbb{P}_n - \mathbb{P})X^T\beta g\big| + \big|2(\mathbb{P}_n - \mathbb{P})fg\big| \\
\triangleq& A + B + C + D + E + F.
\end{aligned} \tag{A.15}$$

We bound each of the terms over the set $\mathcal{M}(R)$ as follows.

*A.* Note that

$$\begin{aligned}
\sup_{(\beta,f,g)\in\mathcal{M}(R)} \big|\|X^T\beta\|_n^2 - \|X^T\beta\|^2\big| &= \sup_{(\beta,f,g)\in\mathcal{M}(R)} \Big|\frac{1}{n}\sum_{i=1}^n \beta^T X_i X_i^T\beta - \beta^T\Sigma\beta\Big| \\
&= \sup_{(\beta,f,g)\in\mathcal{M}(R)} (\beta^T\Sigma\beta)\Big|\frac{1}{n}\sum_{i=1}^n \frac{\beta^T X_i X_i^T\beta}{\beta^T\Sigma\beta} - 1\Big| \\
&\le R_1^2\Big|\frac{1}{n}\sum_{i=1}^n \frac{\beta^T X_i X_i^T\beta}{\beta^T\Sigma\beta} - 1\Big|.
\end{aligned}$$

Note that

$$\mathbb{P}\left(\Big|\frac{1}{n}\sum_{i=1}^{n}\frac{\beta^T X_i X_i^T\beta}{\beta^T\Sigma\beta}-1\Big|>C\left(\sqrt{\frac{t+s_0\log p}{n}}+\frac{t+s_0\log p}{n}\right)\right)\le\exp(-t)$$

for some constant $C>0$. Therefore, by taking $t=n\rho^2/L^2$, we have

$$\sup_{(\beta,f,g)\in\mathcal{M}(R)}\big|\|X^T\beta\|_n^2-\|X^T\beta\|^2\big|\le 2CR_1^2(\rho/L+1/L). \tag{A.16}$$

B. Replace $R^*$ and $K^*$ by $R_1$ and $C_\alpha M_2$, and let $t=n\rho^2/L^2$ in Theorem A.11. Note that

$$\begin{aligned}\mathcal{J}_\infty(C_\alpha M_2,W^{\alpha,2}(R_2,M_2))&\le A_I(R/\rho)^{1/2\alpha}(C_\alpha R/\rho)^{1-1/2\alpha}\\&\le C_\alpha A_I R/\rho\\&\le C_\alpha(\sqrt{n}\rho^{(2\alpha+1)/2\alpha}/L)(R/\rho)\\&\le\sqrt{n}C_\alpha R/L,\end{aligned} \tag{A.17}$$

where the first inequality follows from A.4, the third one and the last one follow from (A.14) and $L>1$.
Then we have with probability at least $1-\exp(-n\rho^2/L)$,

$$\begin{aligned}&\sup_{(\beta,f,g)\in\mathcal{M}(R)}\big|\|f\|_n^2-\|f\|^2\big|\\&\le\frac{2R_1\mathcal{J}_\infty(C_\alpha M_2,W^{\alpha,2}(R_1,M_2))}{\sqrt{n}}+R_1(C_\alpha M_2)\frac{\rho}{L}+\frac{4\mathcal{J}_\infty^2(C_\alpha M_2,W^{\alpha,2}(R_1,M_2))}{n}+(C_\alpha M_2)^2\frac{\rho^2}{L^2}\\&\le(2R_1R/L+R_1R/L+4(R/L)^2+(R/L)^2)C_\alpha^2\\&\le 8C_\alpha^2R_1^2/L,\end{aligned} \tag{A.18}$$

where (A.18) follows from (A.17).

C. Replace $R_1^*$ and $K_1^*$ by $R_1$ and $C_\gamma M_3$, and let $t=n\rho^2/L^2$ in Theorem A.11. By similar arguments as (A.17), together with (A.14), it shows that

$$\mathcal{J}_\infty(C_\gamma M_3,W^{\gamma,2}(R_3,M_3))\le\sqrt{n}C_\gamma R/L. \tag{A.19}$$

Then we have with probability at least $1-\exp(-n\rho^2/L)$,

$$\begin{aligned}&\sup_{(\beta,f,g)\in\mathcal{M}(R)}\big|\|g\|_n^2-\|g\|^2\big|\\&\le\frac{2R_1\mathcal{J}_\infty(C_\gamma M_3,W^{\gamma,2}(R_1,M_3))}{\sqrt{n}}+R_1(C_\gamma M_3)\frac{\rho}{L}+\frac{4\mathcal{J}_\infty^2(C_\gamma M_3,W^{\gamma,2}(R_1,M_3))}{n}+(C_\gamma M_3)^2\frac{\rho^2}{L^2}\\&\le(2R_1R/L+R_1R/L+4(R/L)^2+(R/L)^2)C_\gamma^2\\&\le 8C_\gamma^2R_1^2/L,\end{aligned} \tag{A.20}$$

where (A.20) follow from (A.19) and (A.14).

D. Write $W_i=\beta^TX_i/\|\beta\|$, for $1\le i\le n$. Then we have

$$\sup_{(\beta,f,g)\in\mathcal{M}(R)}\Big|(\mathbb{P}_n-P)(X^T\beta f)\Big|=\sup_{(\beta,f,g)\in\mathcal{M}(R)}\|\beta\|\Big|\frac{1}{n}\sum_{i=1}^{n}W_if(Z_i)\Big|,$$

where $\|\beta\|\le\|\widetilde{X}\beta\|/\Lambda_{\min}\le R/\Lambda_{\min}$ on the set $\mathcal{M}(R)$. Note that $\{W_i,1\le i\le n\}$ are i.i.d. sub-Gaussian with Orlicz norm bounded by $K_X$. Then it follows from Theorem 3.2 of van de Geer (2014)

$$\mathbb{P}\left(\sup_{f\in W^{\alpha,2}(R_1,M_2)}\Big|\frac{1}{n}\sum_{i=1}^{n}W_if(Z_i)\Big|>\frac{2\mathcal{J}_\infty\left(K_XC_\alpha M_2,W^{\alpha,2}(R_1,M_2)\right)+K_XC_\alpha M_2\sqrt{t}}{n}\right)\le\exp(-t). \tag{A.21}$$

By substituting $t = n\rho^2/L^2$, together with (A.17), we have $D \le (2K_X + 2)C_\alpha R^2/(L\Lambda_{\min})$ with probability at least $1 - \exp(-n\rho^2/L^2)$.

E. Similarly as D, we have $E \le \|\beta\| \sup_{(\beta,f,g)\in M(R)} \left|\frac{1}{n}\sum_{i=1}^n W_i g(U_i)\right|$, where

$$\mathbb{P}\left(\sup_{g\in W^{\gamma,2}(R_1,M_3)} \left|\frac{1}{n}\sum_{i=1}^n W_i f(Z_i)\right| > \frac{K_X \mathcal{J}_\infty\left(2C_\gamma M_3, W^{\gamma,2}(R_1,M_3)\right) + K_X C_\gamma M_3 \sqrt{t}}{n}\right) \le \exp(-t).$$

By substituting $t = n\rho^2/L^2$, together with (A.19), we have $E \le (2K_X + 2)C_\gamma R^2/(L\Lambda_{\min})$ with probability at least $1 - \exp(-n\rho^2/L^2)$.

F. Replace $R_1^*, R_2^*, K_1^*, K_2^*$ by $R_1, R_1, C_\alpha M_2, C_\gamma M_3$, and let $t = n\rho^2/L^2$ in Theorem A.12. Note that

$$\begin{aligned} \mathcal{J}_\infty(C_\gamma M_3, W^{\gamma,2}(R_1,M_2)) &\le A_I (R/\rho)^{1/2\gamma}(C_\gamma R/\mu)^{1-1/2\gamma} \\ &\le A_I (R/\rho)^{1/2\gamma}(C_\gamma R/\rho)^{1-1/2\gamma} \\ &\le \sqrt{n} C_\gamma R/L. \end{aligned} \tag{A.22}$$

where the first inequality follows from (A.4), the second one and the last one are from (A.14). Then we obtain with probability at least $1 - \exp(-n\rho^2/L)$,

$$\sup_{(\beta,f,g)\in\mathcal{M}(R)} |\mathbb{P}_n - \mathbb{P})fg| \tag{A.23}$$

$$\begin{aligned} &\le \frac{R_1 \mathcal{J}_\infty(C_\gamma M_3, W^{\gamma,2}(R_1,M_3))}{\sqrt{n}} + \frac{R_1 \mathcal{J}_\infty(R_1(C_\gamma M_3)/R_1, W^{\alpha,2}(R_1,M_2))}{\sqrt{n}} + \frac{R_1(C_\gamma M_3)\rho}{L} \\ &\quad + \frac{C_\alpha M_2 C_\gamma M_3 \rho^2}{L^2} \\ &\le (R_1 R/L + R_1 R/L + R_1 R/L^2 + R^2/L^2) C_\alpha C_\gamma \qquad \text{(A.24)} \\ &\le 4 C_\alpha C_\gamma R_1^2/L, \end{aligned}$$

where(A.24) is implied by (A.22) and (A.14).

Combining $A$ to $F$ and with suitably chosen $L$, we obtain,

$$\sup_{(\beta,f,g)\in\mathcal{M}(R)} \left| \|X^T\beta + f + g\|_n^2 - \|X^T\beta + f + g\|^2 \right| \le \delta_0^2 R^2.$$

with probability at least $1 - 6\exp(-n\rho^2/L)$.

Next, we are going to show with probability at least $1 - 3\exp(-n\rho^2/L)$, the event $\mathcal{T}_2(R)$ hold. Note that $|\mathbb{P}_n\varepsilon(X^T\beta + f + g)| \le \left|\mathbb{P}_n\varepsilon(X^T\beta)\right| + \left|\mathbb{P}_n\varepsilon f\right| + \left|\mathbb{P}_n\varepsilon g\right|$. Note that

$$\sup_{(\beta,f,g)\in\mathcal{M}(R)} \left|\mathbb{P}_n\varepsilon(X^T\beta)\right| \le \sup_{(\beta,f)\in\mathcal{M}(R)} \|\beta\| \left|\frac{1}{n}\sum_{i=1}^n W_i\varepsilon_i\right|. \tag{A.25}$$

Again, we notice that $\{W_i, 1 \le i \le n\}$ are i.i.d. sub-Gaussian and $\{\varepsilon_i, 1 \le i \le n\}$ are i.i.d. sub-Gaussian. It follows from Bernstein's inequality that

$$\mathbb{P}\left(\left|\frac{1}{n}\sum_{i=1}^n W_i\varepsilon_i\right| > C(\sqrt{t/n} + t/n)\right) \le \exp(-t).$$

By taking $t = n\rho^2/L^2$ and together with the fact that $\rho \le R$, we have

$$\sup_{(\beta,f,g)\in\mathcal{M}(R)} \left|\mathbb{P}_n\varepsilon(X^T\beta)\right| \le C(R/\Lambda_{\min})(2\rho/L) \le C(R/\Lambda_{\min})(2R/L),$$

with probability at least $1-\exp(-n\rho^2/L)$. In addition, it follows from Theorem 5.2 of van de Geer and Muro (2015), A.4 and (A.14) that

$$\sup_{(\beta,f,g)\in\mathcal{M}(R)}\left|\mathbb{P}_n\varepsilon f\right| \leq \frac{K_\varepsilon\mathcal{J}_\infty(R_1,W^{\alpha,2}(R_1,M_2))+K_\varepsilon R_1\sqrt{t}}{\sqrt{n}} \leq \frac{K_\varepsilon A_I R}{\sqrt{n}(1-\gamma_1)^{(1-1/2\alpha)/2}\rho^{(1/2\alpha)}}+\frac{R^2}{\sqrt{1-\gamma_1}L}$$
$$\leq \frac{R^2}{L(1-\gamma_1)^{(1-1/2\alpha)/2}}+\frac{R^2}{L\sqrt{1-\gamma_1}} \leq \frac{2R^2}{L\sqrt{1-\gamma_1}},$$

and

$$\sup_{(\beta,f,g)\in\mathcal{M}(R)}\left|\mathbb{P}_n\varepsilon g\right| \leq \frac{K_\varepsilon\mathcal{J}_\infty(R_1,W^{\gamma,2}(R_1,M_3))+K_\varepsilon R_1\sqrt{t}}{\sqrt{n}} \leq \frac{K_\varepsilon A_J M_3^\gamma R^{1-1/2\gamma}}{\sqrt{n}(1-\gamma_1)^{(1-1/2\gamma)/2}}+\frac{R^2}{\sqrt{1-\gamma_1}L}$$
$$\leq \frac{R^2}{L(1-\gamma_1)^{(1-1/2\gamma)/2}}+\frac{R^2}{L\sqrt{1-\gamma_1}} \leq \frac{2R^2}{L\sqrt{1-\gamma_1}}.$$

Therefore, with a suitably chosen $L$, we have

$$\sup_{(\beta,f,g)\in\mathcal{M}(R)}|\mathbb{P}_n\varepsilon(X^T\beta+f+g)| \leq \delta_0^2R^2,$$

with probability at least $1-3\exp(-n\rho^2/L)$. Recalling the probability of $\mathcal{T}_1(R)$, we have shown that for some constants $\widetilde{C}>0, \widetilde{c}>0$, $\mathbb{P}(\mathcal{T}(R)) \geq 1-\widetilde{C}\exp(-\widetilde{c}n\rho^2)$. □

**Lemma A.5.** *Under the conditions of Lemma 2.3, it holds that on* $\mathcal{T}(R)\cap\mathcal{T}_I(R_I)$, $\tau_I(\widehat{\beta}-\beta_0,\widehat{f}-f_0;R_I) \leq R_I$.

***Proof.*** Under the conditions of Lemma 2.3, we can find some $\rho$ and $\mu$ such that

$$\rho^2J_\alpha^2(f_0)+\mu^2J_\gamma^2(g_0) \leq \delta_0^2R^2,\ \rho^2J_\alpha^2(f_0) \leq \delta_I^2R_I^2, \tag{A.26}$$

$$2\mu^2(\Gamma+L_2\delta_0R_I/\lambda)(2\delta_0R/\mu) \leq \delta_IR_I^2,\ \mu^2(\Gamma+L_2\delta_0R_I/\lambda)^2 \leq \delta_I^2 \tag{A.27}$$

for some $\delta_0,\delta_I>0$, which will be taken small enough later.

By the definition of $(\widehat{\beta},\widehat{f},\widehat{g})$, we have

$$\begin{aligned}&\|Y-X^T\widehat{\beta}-\widehat{f}-\widehat{g}\|_n^2+\lambda\|\widehat{\beta}\|_1+\rho^2I^2(\widehat{f})+J_\gamma^2(\widehat{g})\\ \leq&\|Y-X^T\beta_0-f_0-(\widehat{g}+f_{XP}(\widehat{\beta}-\beta_0)+g_X(\widehat{\beta}-\beta_0)+\widehat{f}_P-f_{0P})\|_n^2+\lambda\|\beta_0\|_1+\\ &\rho^2J_\alpha^2(f_0)+\mu^2J_\gamma^2(\widehat{g}+f_{XP}(\widehat{\beta}-\beta_0)+g_X(\widehat{\beta}-\beta_0)+\widehat{f}_P-f_{0P}),\end{aligned}$$

which implies

$$\begin{aligned}&\|\widetilde{X}(\widehat{\beta}-\beta_0)+f_{XA}^T(\widehat{\beta}-\beta_0)+\widehat{f}_A-f_{0A}\|_n^2+\rho^2J_\alpha^2(\widehat{f})\\ \leq&-2\mathbb{P}_n\left(\left((f_{XP}+g_X)^T(\widehat{\beta}-\beta_0)+\widehat{f}_P-f_{0P}+\widehat{g}-g_0\right)\left((\widetilde{X}+f_{XA})^T(\widehat{\beta}-\beta_0)+\widehat{f}_A-f_{0A}\right)\right)\\ &+2\mathbb{P}_n\left(\varepsilon\left(\widetilde{X}+f_{XA}^T(\widehat{\beta}-\beta_0)+\widehat{f}_A-f_{0A}\right)\right)+\lambda\|\beta_0\|_1+\rho^2J_\alpha^2(f_0)-\mu^2J_\gamma^2(\widehat{g})\\ &+\mu^2J_\gamma^2(\widehat{g}+f_{XP}(\widehat{\beta}-\beta_0)+g_X(\widehat{\beta}-\beta_0)+\widehat{f}_P-f_{0P}).\end{aligned}$$

Let

$$t=\frac{R_I}{R_I+\tau_I(\widehat{\beta}-\beta_0,\widehat{f}-f_0;R_I)}.$$

Define $\widetilde{\beta}=t\widehat{\beta}+(1-t)\beta_0$, $\widetilde{f}=t\widehat{f}+(1-t)f_0$, $\widetilde{f}_A=t\widehat{f}_A+(1-t)f_{0A}$. Note that $(\widetilde{\beta},\widetilde{f})\in\mathcal{T}_I(R_I)$ Similarly as the proof of Lemma A.1, it suffices to show that $\tau_I(\widetilde{\beta}-\beta_0,\widetilde{f}-f_0;R_I)\leq R_I/2$.

By convexity and the definition of $\mathcal{T}_I(R_I)$, we have

$$\begin{aligned}&\|\widetilde{X}(\widetilde{\beta}-\beta_0)+f_{XA}^T(\widetilde{\beta}-\beta_0)+\widetilde{f}_A-f_{0A}\|^2+\lambda\|\widetilde{\beta}\|_1+\rho^2J_\alpha^2(\widetilde{f})\\ \le&5R_I^2+\lambda\|\beta_0\|_1+\rho^2J_\alpha^2(f_0)+\mu^2J_\gamma^2(\widetilde{g}+f_{XP}(\widetilde{\beta}-\beta_0)+g_X(\widetilde{\beta}-\beta_0)+\widetilde{f}_P-f_{0P})-\mu^2J_\gamma^2(\widetilde{g}).\end{aligned}$$

Notice that

$$\begin{aligned}&J_\gamma^2(\widetilde{g}+f_{XP}(\widetilde{\beta}-\beta_0)+g_X(\widetilde{\beta}-\beta_0)+\widetilde{f}_P-f_{0P})-J_\gamma^2(\widetilde{g})\\ =&2J_\gamma(\widetilde{g})J_\gamma(f_{XP}^T(\widetilde{\beta}-\beta_0)+g_X^T(\widetilde{\beta}-\beta_0)+\widetilde{f}_P-f_{0P})+J_\gamma^2(f_{XP}^T(\widetilde{\beta}-\beta_0)+g_X^T(\widetilde{\beta}-\beta_0)+\widetilde{f}_P-f_{0P})\\ \le&2J_\gamma(\widetilde{g})\left(J_\gamma(g_X^T(\widetilde{\beta}-\beta_0))+J_\gamma(f_{XP}^T(\widetilde{\beta}-\beta_0)+\widetilde{f}_P-f_{0P})\right)\\ &+\left(J_\gamma(g_X^T(\widetilde{\beta}-\beta_0))+J_\gamma(f_{XP}^T(\widetilde{\beta}-\beta_0)+\widetilde{f}_P-f_{0P})\right)^2\\ \le&2J_\gamma(\widetilde{g})(\|J_\gamma(g_X)\|_\infty\|\widetilde{\beta}-\beta_0\|_1+\Gamma\|f_X(\widetilde{\beta}-\beta_0)+\widetilde{f}-f_0\|)\\ &+(\|J_\gamma(g_X)\|_\infty\|\widetilde{\beta}-\beta_0\|_1+\Gamma\|f_X(\widetilde{\beta}-\beta_0)+\widetilde{f}-f_0\|)^2\\ \le&2\left(\frac{2\delta_0R}{\mu}\right)\left(L_2\frac{\delta_0R_I^2}{\lambda}+\Gamma R_I\right)+\left(L_2\frac{\delta_0R_I^2}{\lambda}+\Gamma R_I\right)^2\\ \le&2\delta_I^2R_I^2/\mu^2.\end{aligned}$$

where the fourth inequality follows from $J(\widetilde{g})\le(2\delta_0R/\mu)$ on $\mathcal{T}(R)$, Assumptions A.4, A.5 and the fact $\|f_X(\widetilde{\beta}-\beta_0)+\widetilde{f}-f_0\|\le R_I$ on $\mathcal{T}_I(R_I)$. The last step follows from (A.27).

Hence, we have

$$\|\widetilde{X}^T(\widetilde{\beta}-\beta_0)\|^2+\|f_{XA}^T(\widetilde{\beta}-\beta_0)+\widetilde{f}_A-f_{0A}\|^2+\lambda\|\widetilde{\beta}\|_1+\rho^2J_\alpha^2(\widetilde{f})\le 8\delta_I^2R_I^2+\lambda\|\beta_0\|_1. \tag{A.28}$$

Subtracting $\lambda\|\widetilde{\beta}_{S_0}\|_1$ on both sides of (A.28), we get

$$\|\widetilde{X}^T(\widetilde{\beta}-\beta_0)\|^2+\|f_{XA}^T(\widetilde{\beta}-\beta_0)+\widetilde{f}_A-f_{0A}\|^2+\lambda\|\widetilde{\beta}_{S_0^c}\|_1+\rho^2J_\alpha^2(\widetilde{f})\le 8\delta_I^2R_I^2+\lambda\|\widetilde{\beta}_{S_0}-\beta_{0S_0}\|_1, \tag{A.29}$$

where $\lambda\|\widetilde{\beta}_{S_0}-\beta_{0S_0}\|_1\le\lambda\sqrt{s_0}\|\beta_{0S_0}-\widetilde{\beta}_{S_0}\|_1\le\lambda\sqrt{s_0}\|\widetilde{\beta}-\beta_0\|_1\le\lambda\sqrt{s_0}\|\widetilde{X}^T(\widetilde{\beta}-\beta_0)\|/\Lambda_{\min}\le\lambda^2s_0/4\Lambda_{\min}^2+\|\widetilde{X}^T(\widetilde{\beta}-\beta_0)\|^2\le\delta_I^2R_I^2+\|\widetilde{X}^T(\widetilde{\beta}-\beta_0)\|^2$. Therefore,

$$\|f_{XA}^T(\widetilde{\beta}-\beta_0)+\widetilde{f}_A-f_{0A}\|^2+\lambda\|\widetilde{\beta}_{S_0^c}\|_1+\rho^2J_\alpha^2(\widetilde{f})\le 9\delta_I^2R_I^2.$$

It holds

($a'$) $\|f_{XA}^T(\widetilde{\beta}-\beta_0)+\widetilde{f}_A-f_{0A}\|\le 3\delta_IR_I$ which further implies

$$\|f_X^T(\widetilde{\beta}-\beta_0)+\widetilde{f}-f_0\|\le 3\delta_IR_I/\sqrt{(1-\gamma^2)};$$

($b'$) $\rho I(\widetilde{f}-f_0)\le\rho I(\widetilde{f})+\rho I(f_0)\le(3+1)\delta_IR_I\le 4\delta_IR_I$ together with equation (A.26).

Note that by using $\lambda\|\widetilde{\beta}_{S_0}-\beta_{0S_0}\|_1\le\lambda^2s_0/(2\Lambda_{\min}^2)+\|X^T(\widetilde{\beta}-\beta_0)\|^2/2$, we can also obtain

($c'$) $\|\widetilde{X}^T(\widetilde{\beta}-\beta_0)\|\le\sqrt{18}\delta_IR_I$.

Now, adding $\lambda\|\widetilde{\beta_{0S_0}}-\beta_{0S_0}\|_1$ on both sides of (A.29), we get

$$\begin{aligned}&\|X^T(\widetilde{\beta}-\beta_0)\|^2+\|f_{XA}^T(\widetilde{\beta}-\beta_0)+\widetilde{f}_A-f_{0A}\|^2+\lambda\|\widetilde{\beta}-\beta_0\|_1+\rho^2J_\alpha^2(\widetilde{f})\\ \le&8\delta_I^2R_I^2+2\lambda\|\widetilde{\beta}_{S_0}-\beta_{0S_0}\|_1\le 8\delta_I^2R_I^2+\lambda^2s_0/\Lambda_{\min}^2+\|X^T(\widetilde{\beta}-\beta_0)\|^2,\end{aligned}$$

which implies that

($d'$) $\lambda\|\widetilde{\beta}-\beta_0\|_1\le 9\delta_I^2R_I^2$.

Combining $(a') - (d')$ and recalling the form of $\tau_I(\widetilde{\beta}-\beta_0, \widetilde{f}-f_0; R_I)$, we obtain

$$\tau_I(\widetilde{\beta}-\beta_0, \widetilde{f}-f_0; R_I) \le ((\sqrt{18}+16)/\sqrt{1-\gamma^2})\delta_I R_I \le \frac{1}{2} R_I,$$

given that $\delta_I \le \sqrt{1-\gamma^2}/(2(\sqrt{18}+16))$. This completes the proof of the lemma.

□

**Lemma A.6.** *Under the conditions of Lemma 2.3, there exist constants $C_I$ and $c_I$, such that*

$$\mathbb{P}(\mathcal{T}_I(R_I)) \ge 1 - C_I \exp(-c_I n\rho^2).$$

***Proof.*** Note that on the set $\mathcal{T}_I(R_I)$, we have $\|f\|^2 \le (1+\Lambda_{\max}/\Lambda'_{\min})R_I^2$, $I(f) \le R_I/\rho$, $\|\widetilde{X}^T\beta\| \le R_I$ and $\|\beta\|_1 \le \delta_0 R_I^2/\lambda$, where $\Lambda'_{\min}$ is defined as the smallest eigenvalue of $\mathbb{E}(f_X f_X^T)$. Also we have $\|g\|^2 \le (1+\Lambda_{\max}/\Lambda_{\min})R^2/(1-\gamma_0)$ and $J(g) \le R/\mu$. Now, we let $R_1'^2 = R_2'^2 = (1+\Lambda_{\max}/\Lambda'_{\min})R_I^2 \triangleq R_I^2/(1-\gamma_1)$ and $R_3'^2 = (1+\Lambda_{\max}/\Lambda_{\min})R^2/(1-\gamma_0) \triangleq R^2/(1-\gamma_2)$, $M_1' = \delta_I R_I^2/\lambda$, $M_2' = R_I/\rho$ and $M_3' = R/\mu$. In particular, we can choose $\rho^2 \le 1-\gamma_1$ and $\mu^2 \le 1-\gamma_2$. Then Lemma A.3 yields that $\sup_{f\in W^{\alpha,2}(R_2',M_2')}\|f\|_\infty \le C_\alpha M_2'$ and $\sup_{g\in W^{\alpha,2}(R_3',M_3')}\|g\|_\infty \le C_\gamma M_3'$. Let $L$ be the constant as in the proof of Lemma A.4. Further, we restrict $\rho, \mu, R$ and $R_I$ as follows

$$R_I \ge RA_J/(\sqrt{n}\mu),\ R_I \ge \rho,\ R_I \ge R\rho/\mu. \tag{A.30}$$

This can be achieved under the assumptions that $\rho^2 \lesssim R_I^2 \le R^2$ and $\mu^2 \lesssim R^2$. We take $t = n\rho^2/L^2$ throughout this proof.

We first look at $\mathcal{T}_{I,1}(R_I)$ and show that with probability at least $1-3\exp(-n\rho^2/L)$, the event $\mathcal{T}_{I,1}(R_I)$ holds. Note that

$$\begin{aligned}
&\|\widetilde{X}^T\beta + f_{XA}^T\beta + f_A\|_n^2 - \|\widetilde{X}^T\beta + f_{XA}^T\beta + f_A\|^2 \\
\le&\left|\|\widetilde{X}^T\beta + f_{XA}^T\beta\|_n^2 - \|\widetilde{X}^T\beta + f_{XA}^T\beta\|^2\right| + \left|\|f_A\|_n^2 - \|f_A\|^2\right| + \left|(\mathbb{P}_n - \mathbb{P})(\widetilde{X}^T\beta + f_{XA}^T\beta)f_A\right| \\
\triangleq& A' + B' + C'.
\end{aligned}$$

We bound $A', B', C'$ as follows, respectively.

$A'$. Recall that $f_{XP}(\cdot) = \mathbb{E}(f_X(Z)|U=\cdot) \in W^{\gamma,2}(L_2)$ and $f_{XA} = f_X - f_{XP}$. We have $f_X, f_{XA}, f_P$ being bounded, therefore sub-Gaussian. Also we have $\|\widetilde{X}^T\beta + f_{XA}^T\beta\|^2 \le 2R_1'^2$. Applying Lemma A.2 and similar arguments as (A.16), we obtain

$$\sup_{(\beta,f)\in\mathcal{M}_I(R_I)} \left|\|\widetilde{X}^T\beta + f_{XA}^T\beta\|_n^2 - \|\widetilde{X}^T\beta + f_{XA}^T\beta\|^2\right| \le 2R_1'^2\left(\frac{\|\widetilde{X}^T\beta + f_{XA}^T\beta\|_n^2}{\|\widetilde{X}^T\beta + f_{XA}^T\beta\|} - 1\right),$$

which is smaller than $4CR_1^{'2}/L$ for some absolute constant $C > 0$ with probability at least $1 - \exp(-n\rho^2/L)$.

$B'$. Note that $\sup_{f\in W^{\alpha,2}(R_2',M_2')}\|f_A\|_\infty \le \sup_{f\in W^{\alpha,2}(R_2',M_2')}\|f_P\|_\infty + \sup_{f\in W^{\alpha,2}(R_2',M_2')}\|f\|_\infty \le 2M_2'$ and $\sup_{f\in W^{\alpha,2}(R_2',M_2')}\|f_A\| \le R_2'$. Replace $R^*$ and $K^*$ by $R_2'$ and $2C_\alpha M_2'$ in Theorem A.11. Note that

$$\mathcal{J}_\infty(2C_\alpha M_2', \{f_A : f\in W^{\alpha,2}(R_2',M_2')\}) \le \mathcal{J}_\infty(2C_\alpha M_2', W^{\alpha,2}(R_2',M_2')) \le 2\sqrt{n}R_I^/L. \tag{A.31}$$

Similarly as (A.18), we then have

$$\begin{aligned}
B' \le& \frac{2R_2'\mathcal{J}_\infty(2C_\alpha M_2', \{f_A : f\in W^{\alpha,2}(R_2',M_2')\})}{\sqrt{n}} + 2R_2'C_\alpha M_2'\frac{\rho}{L} \\
&+ \frac{4\mathcal{J}_\infty^2(2C_\alpha M_2', \{f_A : f\in W^{\alpha,2}(R_2',M_2')\})}{n} + 4C_\alpha^2 M_2^{'2}\frac{\rho^2}{L^2} \\
\le& \left(\frac{8R_I^2}{L\sqrt{(1-\gamma_1)}} + \frac{8R_I^2}{L\sqrt{(1-\gamma_1)}} + \frac{64R_I^2}{L^2} + \frac{4R_I^2}{L}\right)C_\gamma^2 \le \frac{84C_\gamma R_I^2}{L\sqrt{(1-\gamma_1)}},
\end{aligned}$$

where the second inequality follows from (A.31).

$C'$. Write $W_i' = \beta^T(\widetilde{X} + f_{XA})/\|\beta\|, 1 \le i \le n$, which are i.i.d. sub-Gaussian with Orlicz norm bounded by, say $K_X' > 1$.

$$\sup_{(\beta,f)\in\mathcal{M}_I(R_I)} \left|(\mathbb{P}_n - \mathbb{P})(\widetilde{X}^T\beta + f_{XA}^T\beta)f_A\right| = \sup_{(\beta,f)\in\mathcal{M}_I(R_I)} \|\beta\|\left|\frac{1}{n}\sum_{i=1}^n W_i' f_A(Z_i,U_i)\right|,$$

where $\|\beta\| \le \|\widetilde{X}\beta\|/\Lambda_{\min} \le R_I/\Lambda_{\min}$ on the set $\mathcal{M}_I(R_I)$. Similarly as (A.21), we have

$$\sup_{f\in W^{\alpha,2}(R_2',M_2')} \left|\frac{1}{n}\sum_{i=1}^n W_i' f_A(Z_i,U_i)\right| > \frac{2\mathcal{J}_\infty\left(2K_X'C_\gamma M_2', \{f_A : f\in W^{\alpha,2}(R_2',M_2')\}\right) + 2K_X'C_\gamma M_2'\sqrt{t}}{\sqrt{n}}$$

holds with probability $\le \exp(-t)$. Substituting $t = n\rho^2/L$, together with (A.31), we have

$$\sup_{(\beta,f)\in\mathcal{M}_I(R_I)} \left|(\mathbb{P}_n - \mathbb{P})(\widetilde{X}^T\beta + f_{XA}^T\beta)f_A\right| \le C_\gamma(R_I/\Lambda_{\min})(6K_X'R_I/L). \tag{A.32}$$

Combining $A'$ to $C'$, with $L$ large enough, we can have

$$\sup_{(\beta,f)\in\mathcal{M}_I(R_I)} \left|\|\widetilde{X}^T\beta + f_{XA}^T\beta + f_A\|_n^2 - \|\widetilde{X}^T\beta + f_{XA}^T\beta + f_A\|^2\right| \le \delta_I^2 R_I^2$$

with probability at least $1 - 3\exp(-n\rho^2/L)$.

Next, we show with probability at least $1 - 2\exp(-n\rho^2/L)$, the event $\mathcal{T}_{I,2}(R_I)$ holds. Notice that $\left|\mathbb{P}_n(\varepsilon(\widetilde{X}^T\beta + f_{XA}^T\beta + f_A))\right| \le \left|\mathbb{P}_n\varepsilon(\widetilde{X}^T\beta + f_{XA}^T\beta)\right| + \left|\mathbb{P}_n\varepsilon f_A\right|$, where for some absolute constant $C > 0$,

$$\sup_{(\beta,f)\in\mathcal{M}_I(R_I)} \left|\mathbb{P}_n\varepsilon(\widetilde{X}^T\beta + f_{XA}^T\beta)\right| = \sup_{(\beta,f)\in\mathcal{M}_I(R_I)} \|\beta\|\left|\frac{1}{n}\sum_{i=1}^n W_i'\varepsilon_i\right| \le 2CR_I^2/(\Lambda_{\min}L),$$

follows from similar arguments as (A.25). Further, Theorem 5.2 of van de Geer and Muro (2015), A.4 and equation (A.30) shows

$$\begin{aligned}\left|\mathbb{P}_n\varepsilon f_A\right| \le & \frac{K_\varepsilon\mathcal{J}_\infty(R_2', \{f_A : f\in W^{\alpha,2}(R_2',M_2')\}) + K_\varepsilon R_2'\sqrt{t}}{\sqrt{n}} \le \frac{2K_\varepsilon\mathcal{J}_\infty(R_2', W^{\alpha,2}(R_2',M_2')) + K_\varepsilon R_2'\sqrt{t}}{\sqrt{n}} \\ \le & \frac{2K_\varepsilon A_I R_I}{\sqrt{n}(1-\gamma_1)^{(1-1/2\alpha)/2}\rho^{(1/2\alpha)}} + \frac{R_I^2}{\sqrt{1-\gamma_1}L} \le \frac{2R_I^2}{L(1-\gamma_1)^{(1-1/2\alpha)/2}} + \frac{R_I^2}{\sqrt{1-\gamma_1}L} \le \frac{3R_I^2}{\sqrt{1-\gamma_1}L}.\end{aligned}$$

Thus, we have for some suitably chosen $L > 0$,

$$\sup_{(\beta,f)\in\mathcal{M}_I(R_I)} \left|\mathbb{P}_n\varepsilon(\widetilde{X}^T\beta + f_{XA}^T\beta + f_A)\right| \le \delta_I^2 R_I^2$$

with probability at least $1 - 2\exp(-n\rho^2/L)$.

Finally, we show with probability at least $1 - 4\exp(-n\rho^2/L)$, the event $\mathcal{T}_{I,3}(R_I)$ holds. Notice that $\mathbb{E}(\widetilde{X}^T\beta + f_{XA}^T\beta + f_A)(f_{XP}^T\beta + g_X^T\beta + f_P + g) = 0$. Then we get

$$\begin{aligned}&\left|\mathbb{P}_n(\widetilde{X}^T\beta + f_{XA}^T\beta + f_A)(f_{XP}^T\beta + g_X^T\beta + f_P + g)\right| \\ \le&\left|(\mathbb{P}_n - \mathbb{P})\left((\widetilde{X}^T\beta + f_{XA}^T\beta)(f_{XP}^T\beta + g_X^T\beta)\right)\right| + \left|(\mathbb{P}_n - \mathbb{P})(\widetilde{X}^T\beta + f_{XA}^T\beta)(f_P + g)\right| \\ &\left|(\mathbb{P}_n - \mathbb{P})\left((f_{XP}^T\beta + g_X^T\beta)f_A\right)\right| + \left|(\mathbb{P}_n - \mathbb{P})\left(f_A(f_P + g)\right)\right| \\ \triangleq& A'' + B'' + C'' + D''.\end{aligned}$$

It is noted that $\|\widetilde{X}^T\beta + f_{XA}^T\beta\| \le R_1'$, $\|\widetilde{X}^T\beta + f_{XA}^T\beta\|_\infty \le 2M_1'$, $\|f_{XP}^T\beta + g_X^T\beta\| \le R_1'$, $\|f_P + g\| \le \|f_P\| + \|g\| \le R_2' + R_3' \le 2R_3'$, $J(f_P + g) \le J(f_P) + J(g) \le \Gamma\|f\| + (R/\mu) \le \Gamma R_I + (R/\mu) \le 4M_3'$. Then we apply Theorem A.12 for $A'', B'', C'', D''$, respectively. Each of the following terms holds with probability at least $1 - \exp(-n\rho^2/L)$.

$A''$. Note that

$$\begin{aligned}
&\sup_{(\beta,f,g)\in\mathcal{M}(R),(\beta,f)\in\mathcal{M}_I(R_I)} \left|(\mathbb{P}_n-\mathbb{P})\left((\widetilde{X}^T\beta+f_{XA}^T\beta)(f_{XP}^T\beta+g_X^T\beta)\right)\right| \\
&\le \sup_{(\beta,f,g)\in\mathcal{M}(R),(\beta,f)\in\mathcal{M}_I(R_I)} \|\beta\|^2\left|(\mathbb{P}_n-\mathbb{P})\left(\frac{\beta^T(\widetilde{X}+f_{XA})}{\|\beta\|}\frac{\beta^T(f_{XP}+g_X)}{\|\beta\|}\right)\right|,
\end{aligned}$$

where both $\frac{\beta^T(\widetilde{X}+f_{XA})}{\|\beta\|}$ and $\frac{\beta^T(f_{XP}+g_X)}{\|\beta\|}$ are sub-Gaussian. Therefore, by Bernstein's inequality, we have for some constant $C>0$.

$$\mathbb{P}\left(\left|(\mathbb{P}_n-\mathbb{P})\left(\frac{\beta^T(\widetilde{X}+f_{XA})}{\|\beta\|}\frac{\beta^T(f_{XP}+g_X)}{\|\beta\|}\right)\right|>C(\sqrt{t/n}+t/n)\right)\le\exp(-t).$$

By taking $t=n\rho^2/L$ and recalling that $\sup_{(\beta,f)\in\mathcal{M}_I(R_I)}\|\beta\|^2\le R_I^2/\Lambda_{\min}^2$, we can have with probability $1-\exp(-n\rho^2/L)$, that

$$\sup_{(\beta,f,g)\in\mathcal{M}(R),(\beta,f)\in\mathcal{M}_I(R_I)} \left|(\mathbb{P}_n-\mathbb{P})\left((\widetilde{X}^T\beta+f_{XA}^T\beta)(f_{XP}^T\beta+g_X^T\beta)\right)\right|\le 2CR_I^2/(L\Lambda_{\min}^2).$$

$B''$. Recall the definition of $W'$. We have

$$\begin{aligned}
&\sup_{(\beta,f,g)\in\mathcal{M}(R),(\beta,f)\in\mathcal{M}_I(R_I)} \left|(\mathbb{P}_n-\mathbb{P})(\widetilde{X}^T\beta+f_{XA}^T\beta)(f_P+g)\right| \\
&= \sup_{(\beta,f,g)\in\mathcal{M}(R),(\beta,f)\in\mathcal{M}_I(R_I)} \|\beta\|\left|(\mathbb{P}_n-\mathbb{P})W'(f_P+g)\right|.
\end{aligned}$$

Note that on the set $\mathcal{M}_I(R_I)\cap\mathcal{M}(R)$, we have $f_P+g\in W^{\gamma,2}(2R_3',4M_3')$. Further it follows from A.4 and (A.30) that

$$\mathcal{J}_\infty(4C_\gamma M_3', W^{\gamma,2}(2R_3',4M_3'))\le A_JR/\mu\le\sqrt{n}C_\gamma R_I/L$$

and $M_3'\rho/L\le R_I/L$. Therefore, similarly as the proof of (A.32), we can have with probability at least $1-\exp(-n\rho^2/L)$, that

$$\sup_{(\beta,f,g)\in\mathcal{M}(R),(\beta,f)\in\mathcal{M}_I(R_I)} \left|(\mathbb{P}_n-\mathbb{P})(\widetilde{X}^T\beta+f_{XA}^T\beta)(f_P+g)\right|\le C_\gamma(R_I/\Lambda_{\min})(12K_X'R_I/L)$$

$C''$. Write $W''=(f_{XP}^T\beta+g_X^T\beta)/\|\beta\|$, which is sub-Gaussian with Orlicz norm bounded by $K_X''$. Now we have

$$\left|(\mathbb{P}_n-\mathbb{P})\left((f_{XP}^T\beta+g_X^T\beta)f_A\right)\right|=\|\beta\|\left|(\mathbb{P}_n-\mathbb{P})W''f_A\right|.$$

Similarly to the proof of ((A.21)), (A.32), we can have with probability at least $1-\exp(-n\rho^2/L)$,

$$\sup_{(\beta,f)\in\mathcal{M}_I(R_I),(\beta,f,g)\in\mathcal{M}(R)} \|\beta\|\left|(\mathbb{P}_n-\mathbb{P})W''f_A\right|\le C_\alpha(R_I/\Lambda_{\min})(10K_X''R_I/L).$$

$D''$. Similar to the proof of (A.23), we have

$$\begin{aligned}
D''&\le\frac{R_2'\mathcal{J}_\infty(4C_\gamma M_3',W^{\gamma,2}(2R_3',4M_3'))}{\sqrt{n}}+\frac{2R_3'\mathcal{J}_\infty(R_2'(4C_\gamma M_3')/2R_3',\{f_A:f\in W^{\alpha,2}(R_2',M_2')\})}{\sqrt{n}} \\
&\quad+\frac{R_2'(4C_\gamma M_3')\rho}{L}+\frac{(2C_\alpha M_2')(4C_\gamma M_3')\rho^2}{L^2} \\
&\le\Big(\frac{R_I}{\sqrt{1-\gamma_1}}\frac{4R_I}{L}+\left(\frac{\sqrt{1-\gamma_2}}{\sqrt{1-\gamma_1}}\right)^{1-1/2\alpha}\frac{8R^{(1/2\alpha)}}{\sqrt{1-\gamma_2}}\frac{A_I(R/\mu)^{1-1/2\alpha}R_I\rho}{\sqrt{n}\rho^{1+1/2\alpha}}+\frac{R_I}{\sqrt{1-\gamma_1}}\frac{4(R/\mu)\rho}{L} \\
&\quad+\frac{2R_I}{\rho}\frac{4(R/\mu)\rho^2}{L^2}\Big)C_\alpha C_\gamma \\
&\le\frac{18C_\alpha C_\gamma R_I^2}{\sqrt{1-\gamma_1}(\sqrt{1-\gamma_2})^{1/2\alpha}L}.
\end{aligned}$$

Therefore, by choosing $L$ large enough, we can have with probability at least $1 - 4\exp(-n\rho^2/L)$,

$$\sup_{(\beta,f,g)\in\mathcal{M}(R),(\beta,f)\in\mathcal{M}_I(R_I)} \left|\mathbb{P}_n(\widetilde{X}^T\beta + f_{XA}^T\beta + f_A)(f_{XP}^T\beta + g_X^T\beta + f_P + g)\right| \le \delta_I^2 R_I^2,$$

by letting $L$ large enough. Now, we conclude that there exists constant $C_I$ and $c_I$, such that

$$\mathbb{P}(\mathcal{T}_I(R_I)) \ge 1 - C_I\exp(-c_I n\rho^2).$$

□

***Proof of Lemma 2.3.*** This proof simply follows from the following inequality

$$\mathbb{P}(\mathcal{T}(R)\cap\mathcal{T}_I(R_I)) \ge 1 - \mathbb{P}(\mathcal{T}^c(R)) - \mathbb{P}(\mathcal{T}_I^c(R_I))$$

and Lemma A.4 and A.6. □

*A.2.2. Proof of Lemma 2.4*

We start from the main proof of Lemma 2.4, followed by some necessary lemmas.

***Proof.*** Recall that $\pi_{X|Z,U} = f_X + g_X$. By the definition of $(\widehat{\beta},\widehat{f},\widehat{g})$, we have

$$\begin{aligned}
&\|Y - X^T\widehat{\beta} - \widehat{f} - \widehat{g}\|_n^2 + \lambda\|\widehat{\beta}\|_1 + \rho^2 J_\alpha^2(\widehat{f}) + \mu^2 J_\gamma^2(\widehat{g})\\
\le&\|Y - X^T\beta_0 - (\widehat{f} + f_X^T(\widehat{\beta}-\beta_0)) - (\widehat{g} + g_X^T(\widehat{\beta}-\beta_0))\|_n^2 + \lambda\|\beta_0\|_1 + \rho^2 J_\alpha^2(\widehat{f} + f_X^T(\widehat{\beta}-\beta_0))\\
&+ \mu^2 J_\gamma^2(\widehat{g} + g_X^T(\widehat{\beta}-\beta_0)).
\end{aligned}$$

That implies

$$\begin{aligned}
&\|\widetilde{X}(\widehat{\beta}-\beta_0)\|_n^2 + \lambda\|\widehat{\beta}\|_1 + \rho^2 J_\alpha^2(\widehat{f}) + \mu^2 J_\gamma^2(\widehat{g}) \qquad (A.33)\\
\le&2\left|\mathbb{P}_n\varepsilon(\widetilde{X}^T(\widehat{\beta}-\beta_0))\right| + 2\left|\mathbb{P}_n(\widehat{f} - f_0 + \widehat{g} - g_0)\widetilde{X}^T(\widehat{\beta}-\beta_0)\right|\\
&+ 2\left|\mathbb{P}_n(f_X^T(\widehat{\beta}-\beta_0) + g_X^T(\widehat{\beta}-\beta_0))\widetilde{X}^T(\widehat{\beta}-\beta_0)\right| + \lambda\|\beta_0\|_1 + \rho^2 J_\alpha^2(\widehat{f} + f_X^T(\widehat{\beta}-\beta_0))\\
&+ \mu^2 J_\gamma^2(\widehat{g} + g_X^T(\widehat{\beta}-\beta_0)).
\end{aligned}$$

From Lemmas A.7-A.10, we know that with probability $1 - 7/(2p) - c\exp(-Cn\rho^2)$ for some constant $c, C > 0$, (A.33) can be further reduced to

$$\|\widetilde{X}(\widehat{\beta}-\beta_0)\|_n^2 + \lambda\|\widehat{\beta}\|_1 \le (\lambda/2)\|\widehat{\beta}-\beta_0\|_1 + \lambda\|\beta_0\|_1.$$

Hence, Noting that $\|\widehat{\beta}\|_1 = \|\widehat{\beta}_{S_0}\|_1 + \|\widehat{\beta}_{S_0^c}\|_1$ and $\|\beta_0\|_1 = \|\beta_{0S_0}\|_1$, we get $\|\widetilde{X}(\widehat{\beta}-\beta_0)\|^2 + \frac{1}{2}\lambda\|\widehat{\beta}_{S_0^c}\|_1 \le \frac{\lambda}{2}\|\widehat{\beta}_{S_0} - \beta_{0S_0}\|_1 + \frac{\lambda}{2}\|\widehat{\beta}_{S_0^c} - \beta_{0S_0^c}\|_1 + \lambda\|\beta_{0S_0}\|_1 - \lambda\|\widehat{\beta}_{S_0}\|_1 - \frac{\lambda}{2}\|\widehat{\beta}_{S_0^c}\|_1 \le \frac{3\lambda}{2}\|\widehat{\beta}_{S_0} - \beta_{0S_0}\|_1$. This gives

$$\begin{aligned}
&2\|\widetilde{X}(\widehat{\beta}-\beta_0)\|_n^2 + \lambda\|\widehat{\beta}-\beta_0\|_1 = 2\|\widetilde{X}(\widehat{\beta}-\beta_0)\|_n^2 + \lambda\|\widehat{\beta}_{S_0^c}\|_1 + \lambda\|\widehat{\beta}_{S_0^c} - \beta_{0S_0^c}\|_1\\
\le&4\lambda\|\widehat{\beta}_{S_0} - \beta_{0S_0}\|_1 \le 4\lambda\sqrt{s_0}\|\widehat{\beta}-\beta_0\| \le \frac{4\lambda\sqrt{s_0}}{\Lambda_{\min}}\|\widetilde{X}(\widehat{\beta}-\beta_0)\| \le \|\widetilde{X}(\widehat{\beta}-\beta_0)\|^2 + \frac{4\lambda^2 s_0}{\Lambda_{\min}^2}\\
\le&\|\widetilde{X}(\widehat{\beta}-\beta_0)\|_n^2 + \frac{4\lambda^2 s_0}{\Lambda_{\min}^2} + \frac{\lambda}{2}\|\widehat{\beta}-\beta_0\|_1.
\end{aligned}$$

Therefore,

$$\|\widetilde{X}(\widehat{\beta}-\beta_0)\|_n^2 + \frac{\lambda}{2}\|\widehat{\beta}-\beta_0\|_1 \le \frac{4\lambda^2 s_0}{\Lambda_{\min}^2}.$$

□

**Lemma A.7.** *With probability at least* $1-1/p$,

$$2\left|\mathbb{P}_n\varepsilon(\widetilde{X}^T(\widehat{\beta}-\beta_0))\right| \leq 4\sqrt{6K_{\widetilde{X}}}K_\varepsilon\sqrt{\frac{\log(2p)}{n}}\|\widehat{\beta}-\beta_0\|_1 \leq \frac{\lambda}{10}\|\widehat{\beta}-\beta_0\|_1$$

*for some constant* $K_{\widetilde{X}} > 1$.

***Proof.*** First we have

$$\left|\mathbb{P}_n\varepsilon(\widetilde{X}^T(\widehat{\beta}-\beta_0))\right| \leq \|\mathbb{P}_n\varepsilon\widetilde{X}^T\|_\infty\|\widehat{\beta}-\beta_0\|_1.$$

Assumption A.2 that $\mathbb{E}\exp(\varepsilon_i^2/K_\varepsilon^2) \leq 2$ implies $\mathbb{E}\exp(t\varepsilon_i) \leq \exp(3K_\varepsilon^2 t/2)$, see Vershynin (2012). Then we get

$$\mathbb{E}\exp(t(\frac{1}{n}\sum_{i=1}^n\varepsilon_i\widetilde{X}_{ij})) = \prod_{i=1}^n\mathbb{E}\exp(\frac{t}{n}\widetilde{X}_{ij}\varepsilon_i) \leq \prod_{i=1}^n\exp(\frac{3}{2}K_\varepsilon^2(\frac{t^2}{n^2}\widetilde{X}_{ij}^2)) = \exp(\frac{3}{2}K_\varepsilon^2\frac{t^2}{n}\|\widetilde{X}_j\|_n^2),$$

which implies given $\widetilde{X}$ fixed, for $t>0$ and all $j$,

$$\mathbb{P}\left\{\left|\frac{1}{n}\sum_{i=1}^n\varepsilon_i\widetilde{X}_{ij}\right| > \sqrt{\frac{t}{n}}2\|\widetilde{X}_j\|_n\sqrt{\frac{3}{2}}K_\varepsilon\right\} \leq \exp(-t),$$

see Vershynin (2012). Hence

$$\mathbb{P}\left\{\max_{1\leq j\leq p}\left|\frac{1}{n}\sum_{i=1}^n\varepsilon_i\widetilde{X}_{ij}\right| > \sqrt{\frac{t+\log p}{n}}2\|\widetilde{X}_j\|_n\sqrt{\frac{3}{2}}K_\varepsilon\right\} \leq \exp(-t).$$

Note that $\pi_{X|Z,U} = f_X(Z) + g_X(U)$ with $f_X \in W^{\alpha,2}(L_1)$ and $g_X \in W^{\gamma,2}(L_2)$. Therefore we have $\widetilde{X} = X - \pi_{X|Z,U}$ is sub-Gaussian. Then by Bernstein's inequality, we have for some constant $K_{\widetilde{X}} \geq 1$ that

$$\mathbb{P}\left\{\max_{1\leq j\leq p}\left|\|\widetilde{X}_j\|_n^2 - \mathbb{E}\|\widetilde{X}_j\|_n^2\right| \geq K_{\widetilde{X}}\sqrt{\frac{\log(2p)}{n}}\right\} \leq 1/(2p),$$

which further implies

$$\mathbb{P}\left\{\max_{1\leq j\leq p}\|\widetilde{X}_j\|_n^2 \geq 2K_{\widetilde{X}}\right\} \leq \mathbb{P}\left\{\max_{1\leq j\leq p}\|\widetilde{X}_j\|_n^2 \geq \mathbb{E}\|\widetilde{X}_j\|_n^2 + K_{\widetilde{X}}\sqrt{\frac{\log(2p)}{n}}\right\} \leq 1/(2p), \tag{A.34}$$

Now take $t = \log(2p)$. With probability at least $1-1/p$,

$$\|\mathbb{P}_n\varepsilon\widetilde{X}^T\|_\infty \leq \sqrt{6}\sqrt{2K_{\widetilde{X}}}K_\varepsilon\sqrt{\frac{2\log(2p)}{n}}.$$

Noting that $\lambda \gtrsim \sqrt{\log p/n}$, we can have $2\|\mathbb{P}_n\varepsilon\widetilde{X}^T\|_\infty \leq 4\sqrt{6K_{\widetilde{X}}}K_\varepsilon\sqrt{\log(2p)/n} \leq \lambda/10$.

□

**Lemma A.8.** *With probability at least* $1 - 5/(2p) - C\exp(-n\rho^2/c)$ *for some constants* $c, C > 0$,

$$2\left|\mathbb{P}_n(\widehat{f}-f_0+\widehat{g}-g_0)\widetilde{X}^T(\widehat{\beta}-\beta_0)\right| + 2\left|\mathbb{P}_n(f_X^T(\widehat{\beta}-\beta_0)+g_X^T(\widehat{\beta}-\beta_0))\widetilde{X}^T(\widehat{\beta}-\beta_0)\right| \leq \frac{\lambda}{10}\|\widehat{\beta}-\beta_0\|_1,$$

$$\left|\|\widetilde{X}(\widehat{\beta}-\beta_0)\|_n^2 - \|\widetilde{X}(\widehat{\beta}-\beta_0)\|^2\right| \leq \frac{\lambda}{2}\|\widehat{\beta}-\beta_0\|_1.$$

***Proof.*** On the set $\mathcal{T}(R)$, we have

$$\begin{aligned}\left|\mathbb{P}_n(\widehat{f}-f_0+\widehat{g}-g_0)\widetilde{X}^T(\widehat{\beta}-\beta_0)\right| &\leq \left(\|\mathbb{P}_n(\widehat{f}-f_0+\widehat{g}-g_0)\widetilde{X}^T\|_\infty\right)\|\widehat{\beta}-\beta_0\|_1\\ &\leq \max_{1\leq j\leq p}\left(\Big|\frac{1}{n}\sum_{i=1}^n(\widehat{f}-f_0+\widehat{g}-g_0)_i\widetilde{X}_{ij}\Big|\right)\|\widehat{\beta}-\beta_0\|_1.\end{aligned}$$

Note that given $\widetilde{X}$, we have for each $1\leq j\leq p$,

$$|(\widehat{f}-f_0+\widehat{g}-g_0)_i\widetilde{X}_{ij}|\leq 2R/\sqrt{1-\gamma_2}|\widetilde{X}_{ij}|.$$

By Lemma 14.15 of Bühlmann and van de Geer (2011), we have

$$\mathbb{P}\left\{\max_{1\leq j\leq p}\Big|\frac{1}{n}\sum_{i=1}^n(\widehat{f}-f_0+\widehat{g}-g_0)_i\widetilde{X}_{ij}\Big|\geq\max_{1\leq j\leq p}\sqrt{\frac{(2R/\sqrt{1-\gamma_2})^2\sum_{i=1}^n\widetilde{X}_{ij}^2}{n}}\sqrt{2\left(t^2+\frac{\log(p)}{n}\right)}\right\}$$
$$\leq\exp(-nt^2).$$

Take $t^2=\delta_3^2\rho^2$ for some $\delta_3>0$. Again by (A.34) and noticing that $\delta_3\rho>\sqrt{\log p/n}$, we have that with probability at least $1-1/(2p)-\exp(-n\delta_3\rho^2)$,

$$\max_{1\leq j\leq p}\Big|\frac{1}{n}\sum_{i=1}^n(\widehat{f}-f_0+\widehat{g}-g_0)_i\widetilde{X}_{ij}\Big|\leq\sqrt{2K_X}(2R/\sqrt{1-\gamma_2})2\delta_3\mu.$$

Therefore, by choosing $\delta_3$ suitably small, we can have

$$\max_{1\leq j\leq p}\Big|\frac{1}{n}\sum_{i=1}^n(\widehat{f}-f_0+\widehat{g}-g_0)_i\widetilde{X}_{ij}\Big|\leq\lambda/20.$$

Next, note that

$$\begin{aligned}&\Big|\frac{1}{n}\sum_{i=1}^n((f_X+g_X)^T(\widehat{\beta}-\beta_0))_i\widetilde{X}_{ij}\Big|\leq\Big|\frac{1}{n}\sum_{i=1}^n\sum_{k=1}^p(f_X+g_X)_{ik}(\widehat{\beta}-\beta_0)_k\widetilde{X}_{ij}\Big|\\ \leq&\Big|\sum_{k=1}^p(\widehat{\beta}-\beta_0)_k\Big(\frac{1}{n}\sum_{i=1}^n(f_X+g_X)_{ik}\widetilde{X}_{ij}\Big)\Big|\leq\|\widehat{\beta}-\beta_0\|_1\max_{1\leq k\leq p}\Big|\frac{1}{n}\sum_{i=1}^n(f_X+g_X)_{ik}\widetilde{X}_{ij}\Big|\\ \leq&\delta_0\frac{R^2}{\lambda}\max_{1\leq k\leq p}\Big|\frac{1}{n}\sum_{i=1}^n(f_X+g_X)_{ik}\widetilde{X}_{ij}\Big|,\end{aligned}$$

where $\mathbb{E}(f_X+g_X)_{ik}\widetilde{X}_{ij}=0$ and $|(f_X+g_X)_{ik}\widetilde{X}_{ij}|\leq M_0|\widetilde{X}_{ij}|$ conditional on $\widetilde{X}$. By Lemma 14.15 in Bühlmann and van de Geer (2011), we obtain that given $\widetilde{X}$,

$$\mathbb{P}\left(\max_{1\leq j\leq p}\max_{1\leq k\leq p}\Big|\frac{1}{n}\sum_{i=1}^n(f_X+g_X)_{ik}\widetilde{X}_{ij}\Big|\geq\max_{1\leq j\leq p}\sqrt{\frac{M_0^2\sum_{i=1}^n\widetilde{X}_{ij}^2}{n}}\sqrt{2\left(t^2+\frac{2\log 2p}{n}\right)}\right)$$
$$\leq\exp(-nt^2).$$

Similarly, letting $t^2=\log(2p)/n$ and revoking (A.34) gives

$$\mathbb{P}\left(\max_{1\leq j\leq p}\max_{1\leq k\leq p}\Big|\frac{1}{n}\sum_{i=1}^n(f_X+g_X)_{ik}\widetilde{X}_{ij}\Big|>\sqrt{2K_X}M_0\sqrt{\frac{\log 2p}{n}}\right)\leq 1/p.$$

Choose $\lambda > 2K_X M_0\sqrt{\log(2p)/n}$. We finally get with probability at least $1-1/p$,

$$2\left|\mathbb{P}_n(f_X^T(\widehat{\beta}-\beta_0)+g_X^T(\widehat{\beta}-\beta_0))\widetilde{X}^T(\widehat{\beta}-\beta_0)\right| \le \delta_0 R^2\|\widehat{\beta}-\beta_0\|_1$$

which can be smaller than $\frac{\lambda}{20}\|\widehat{\beta}-\beta_0\|_1$ by taking suitable choices of $\delta_0$.

Now, we show the second part of the Lemma. Similarly, we have on the set $\mathcal{T}(R)$,

$$|\|\widetilde{X}(\widehat{\beta}-\beta_0)\|_n^2 - \|\widetilde{X}(\widehat{\beta}-\beta_0)\|^2| \le \delta_0\frac{R^2}{\lambda}\max_{1\le k,j\le p}\left|\frac{1}{n}\sum_{i=1}^n(\widetilde{X}_{ik}\widetilde{X}_{jk}-\mathbb{E}\widetilde{X}_{ik}\widetilde{X}_{jk})\right|\|\widehat{\beta}-\beta_0\|_1,$$

where $\widetilde{X}_{ik}\widetilde{X}_{jk}-\mathbb{E}\widetilde{X}_{ik}\widetilde{X}_{jk}$ is sub-exponential. By Bernstein's inequality, we have for some constant $K_{\widetilde{X}}$ that

$$\mathbb{P}\left(\max_{1\le j,k\le p}\left|\frac{1}{n}\sum_{i=1}^n(\widetilde{X}_{ik}\widetilde{X}_{jk}-\mathbb{E}\widetilde{X}_{ik}\widetilde{X}_{jk})\right| > K_{\widetilde{X}}\sqrt{\frac{\log 2p}{n}}\right) \le 1/(2p).$$

Therefore, by choosing $\lambda > 2\delta_0 K_{\widetilde{X}}\sqrt{\log 2p/n}$, we have $|\|\widetilde{X}(\widehat{\beta}-\beta_0)\|_n^2 - \|\widetilde{X}(\widehat{\beta}-\beta_0)\|^2| \le \lambda\|\widehat{\beta}-\beta_0\|_1/2$, with probability at least $1-1/(2p)$. Recalling the probability of $\mathcal{T}(R)$ from Lemma A.4, this lemma is proved.

□

**Lemma A.9.** *Assume*

$$\rho^2 \le \frac{\delta_0^2 R^2}{2(L_1+J_\alpha(f_0))L_1}.$$

*Then on the set* $\mathcal{T}(R)$,

$$\left|\rho^2 J_\alpha^2(\widehat{f}+f_X^T(\widehat{\beta}-\beta_0)) - \rho^2 J_\alpha^2(\widehat{f})\right| \le \frac{\lambda}{10}\|\widehat{\beta}-\beta_0\|_1.$$

***Proof.***

$$\begin{aligned}
\rho^2 J_\alpha^2(\widehat{f}+f_X^T(\widehat{\beta}-\beta_0)) - \rho^2 J_\alpha^2(\widehat{f}) &= \rho^2[J_\alpha^2(f_X^T(\widehat{\beta}-\beta_0)) + 2J_\alpha(\widehat{f}, f_X^T(\widehat{\beta}-\beta_0)] \\
&\le \rho^2\frac{\delta_0 R^2}{\lambda}L_1^2\|\widehat{\beta}-\beta_0\|_1 + 2J_\alpha(\widehat{f})J_\alpha(f_X^T(\widehat{\beta}-\beta_0) \\
&\le \left[\delta_0\rho^2 L_1^2 + 2\rho^2(\frac{R}{\rho}+J_\alpha(f_0))L_1\right]\|\widehat{\beta}-\beta_0\|_1 \\
&\le \left(\frac{1}{2}\delta_0^3 R^2 + \frac{\sqrt{2}}{2}\delta_0 R^2 + \delta_0^2 R^2\right)\|\widehat{\beta}-\beta_0\|_1 \\
&\le 3\delta_0 R^2\|\widehat{\beta}-\beta_0\|_1,
\end{aligned}$$

where the first equality follows from definition of $J_\alpha(\cdot)$, the second inequality follows from Assumption A.5 and $\|\widehat{\beta}-\beta_0\|_1 \le \frac{\delta_0 R^2}{\lambda}$ on $\mathcal{T}(R)$, and the third one is true due to triangular inequality. Choosing $\delta_0$ such that $3\delta_0 R^2 \le \lambda/10$, we get the desired result. □

**Lemma A.10.** *Assume*

$$\mu^2 \le \frac{\delta_0^2 R^2}{(L_2+J_\gamma(g_0))L_2}. \tag{A.35}$$

*Then on the set* $\mathcal{T}(R)$,

$$\left|\mu^2 J_\gamma^2(\widehat{g}+g_X^T(\widehat{\beta}-\beta_0)) - \mu^2 J_\gamma^2(\widehat{g})\right| \le \frac{\lambda}{10}\|\widehat{\beta}-\beta_0\|_1$$

***Proof.***

$$\mu^2 J_\gamma^2(\widehat{g} + g_X^T(\widehat{\beta} - \beta_0)) - \mu^2 J_\gamma^2(\widehat{g}) \le \mu^2 \left(2J(\widehat{g})J(g_X^T(\widehat{\beta} - \beta_0)) + 2J_\gamma^2(g_X^T(\widehat{\beta} - \beta_0))\right)$$

Note that $J(g_X^T(\widehat{\beta} - \beta_0)) \le L_2\|\widehat{\beta} - \beta_0\|_1$ by Assumption A.5 and $\|\widehat{\beta} - \beta_0\|_1 \le \frac{\delta_0 R^2}{\lambda}$ on $\mathcal{T}(R)$. We have

$$\begin{aligned}
\mu^2 J_\gamma^2(\widehat{g} + g_X^T(\widehat{\beta} - \beta_0)) - \mu^2 J_\gamma^2(\widehat{g}) =& \mu^2[2J_\gamma(\widehat{g}, g_X^T(\widehat{\beta} - \beta_0) + J_\gamma^2(g_X^T(\widehat{\beta} - \beta_0))] \\
\le& \mu^2 \left(2\left(R/\mu + J_\gamma(g_0)\right) L_2\|\widehat{\beta} - \beta_0\|_1 + L_2^2\|\widehat{\beta} - \beta_0\|_1^2\right) \\
\le& \left(2\delta_0 R^2 + 2\delta_0^2 R^2 + \delta_0^2 R^2 \left(\frac{\delta_0 R^2}{\lambda}\right)\right)\|\widehat{\beta} - \beta_0\|_1 \\
\le& 5\delta_0 R^2\|\widehat{\beta} - \beta_0\|_1,
\end{aligned}$$

where the first inequality follows from definition and the second one follows from the condition (A.35). Choosing $\delta_0$ such that $5\delta_0 R^2 \le \lambda/10$, we get the desired result. □

*A.2.3. Proof of Theorem 2.5*

***Proof.*** Let $\lambda = 4\sqrt{\log p/n}$, $\rho^2 = n^{-2\alpha/(2\alpha+1)}$, $\mu^2 = n^{-2\gamma/(2\gamma+1)}$. And recall that $R^2 \asymp \mu^2 + \lambda^2 s_0$ and $R_I^2 \asymp \rho^2 + \lambda^2 s_0$. It follows from the definition of $(\widehat{\beta}, \widehat{f}, \widehat{g})$ that

$$\|Y - X\widehat{\beta} - \widehat{f} - \widehat{g}\|_n^2 + \lambda\|\widehat{\beta}\|_1 + \rho^2 J_\alpha^2(\widehat{f}) + \mu^2 J_\gamma^2(\widehat{g}) \le \|Y - X^T\beta_0 - f_0\|_n^2 + \lambda\|\beta_0\|_1 + \rho^2 J_\alpha^2(f_0) + \mu^2 J_\gamma^2(g_0). \quad \text{(A.36)}$$

We first show the risk bound for $\widehat{\beta}$. Triangle inequality and (A.36) imply

$$\lambda\|\widehat{\beta} - \beta_0\|_1 \le \lambda\|\widehat{\beta}\|_1 + \lambda\|\widehat{\beta}_0\|_1 \le \|\varepsilon\|_n^2 + 2\lambda\|\beta_0\|_1 + \rho^2 J_\alpha^2(f_0) + \mu^2 J_\gamma^2(g_0),$$

which further implies for any $k \ge 1$,

$$\mathbb{E}\|\widehat{\beta} - \beta_0\|^{(1/2\alpha)} \le \mathbb{E}\|\widehat{\beta} - \beta_0\|_1^{(1/2\alpha)} \le \mathbb{E}(\|\varepsilon\|_n^2/\lambda + \|\beta_0\|_1 + \rho^2 J_\alpha^2(f_0)/\lambda + \mu^2 J_\gamma^2(g_0)/\lambda)^{(1/2\alpha)}.$$

Note that $n\|\varepsilon\|_n^2$ follows chi-squared distribution with degree of freedom $n$. Thus we have $\mathbb{E}\|\varepsilon\|_n^{(1/2\alpha)} = O(1)$. Also we have that $\|\beta_0\|_1 = O(\sqrt{s_0})$. Therefore, it follows that

$$\mathbb{E}\|\widehat{\beta} - \beta_0\|^{(1/2\alpha)} \le \mathbb{E}\|\widehat{\beta} - \beta_0\|_1^{(1/2\alpha)} \le O(1/\lambda + \sqrt{s_0})^{(1/2\alpha)}.$$

Define the set $\mathcal{T}_1 = \{\|\widetilde{X}(\widehat{\beta} - \beta_0)\|^2 \le \lambda^2 s_0\}$. Then it's known from the proof of Lemma 2.4 that $\mathbb{P}(\mathcal{T}_1^c) \le c/p^4 + c\exp(-n\rho^2/c)$ for some constant $c > 0$. Hence we have

$$\begin{aligned}
\mathbb{E}\|\widehat{\beta} - \beta_0\|^2 &= \mathbb{E}\|\widehat{\beta} - \beta_0\|^2 1_{\mathcal{T}_1} + \mathbb{E}\|\widehat{\beta} - \beta_0\|^2 1_{\mathcal{T}_1^c} \\
&\le O(\lambda^2 s_0) + \sqrt{\mathbb{E}\|\widehat{\beta} - \beta_0\|^4}\sqrt{\mathbb{P}(\mathcal{T}_1^c)} \\
&\le O(\lambda^2 s_0) + O(1/\lambda^2 + s_0)\sqrt{\exp(-n\rho^2) + 1/p^4} \\
&\le O(\lambda^2 s_0).
\end{aligned}$$

The last inequality holds due to the following arguments,

(i). $(1/\lambda^2)\sqrt{\exp(-n^{1/(2\alpha+1)})} = O(\lambda^2 s_0)$, since $n^2\exp(-n^{1/(2\alpha+1)}) = O(s_0 \log^2 p)$;
(ii). $(1/\lambda^2)(1/p^2) = O(\lambda^2 s_0)$, since $1 = O((p^2 s_0 \log^2 p)/n^2)$;
(iii). $s_0\sqrt{\exp(-n^{1/(2\alpha+1)})} = O(\lambda^2 s_0)$, since $n\exp(-n^{1/(2\alpha+1)}) = O(\log p)$;
(iv). $s_0(1/p^2) = O(\lambda^2 s_0)$, since $1 = O((p^2 \log p)/n)$.

We next show the estimation risk for $\widehat{f}$. Define $\mathcal{T}_2 = \{\tau(\widehat{\beta}-\beta_0, \widehat{f}-f_0, \widehat{g}-g_0; R) \le R, \tau_I(\widehat{\beta}-\beta_0, \widehat{f}-f_0; R_I) \le R_I\}$. Note that $\widehat{f} \in W^{2,\alpha}(L_1)$, then together with (A.36), it implies that, for some constant $C>0$, $\sup_{z\in[0,1]} |\widehat{f}(z)-f_0(z)|^2 \le C J_\alpha^2(\widehat{f}-f_0) = C\int_0^1 (\widehat{f}^{(\alpha)}(z) - f_0^{(\alpha)}(z))^2 dz \le C(\|\varepsilon\|_n^2/\rho^2 + 2\lambda\|\beta_0\|_1/\rho^2 + J_\alpha^2(f_0) + \mu^2 J_\gamma^2(g_0)/\rho^2)$. Therefore we have

$$\mathbb{E}\|\widehat{f}-f_0\|_n^{2k} \le (\mathbb{E}\int_0^1 (\widehat{f}^{(\alpha)}(z) - f_0^{(\alpha)}(z))^2 dz)^{(1/2\alpha)} \le O(1/\rho^2 + \lambda\sqrt{s_0}/\rho^2 + \mu^2/\rho^2)^{(1/2\alpha)} = O(1/\rho^{2k}), \tag{A.37}$$

for any $k \ge 1$. Hence, we have

$$\begin{aligned}\mathbb{E}\|\widehat{f}-f_0\|_n^2 &= \mathbb{E}\|\widehat{f}-f_0\|_n^2 1_{\mathcal{T}_2} + \mathbb{E}\|\widehat{f}-f_0\|_n^2 1_{\mathcal{T}_2^c} \\ &\le O(\rho^2+\lambda^2 s_0) + \sqrt{\mathbb{E}\|\widehat{f}-f_0\|_n^4}\sqrt{\mathbb{P}(\mathcal{T}_2^c)} \\ &\le O(\rho^2+\lambda^2 s_0) + O(1/\rho^2)\sqrt{\exp(-nc\rho^2)} \\ &\le O(\rho^2+\lambda^2 s_0),\end{aligned}$$

where the second inequality follows from Lemma 2.3 and the last step is true since $\sqrt{\exp(-n^{1/(2\alpha+1)})} = O(n^{-4/(2\alpha+1)})$.

Now, we are going to show the risk bound for $\widehat{g}$. Define $\mathcal{T}_3 = \{\tau(\widehat{\beta}-\beta_0, \widehat{f}-f_0, \widehat{g}-g_0; R) \le R\}$. By similar arguments as (A.37), we have $\mathbb{E}\|\widehat{g}-g_0\|_n^{2k} = O(1/\mu^{2k})$. Then together with Lemma 2.3, it shows

$$\begin{aligned}\mathbb{E}\|\widehat{g}-g_0\|_n^2 = \mathbb{E}\|\widehat{g}-g_0\|_n^2 1_{\mathcal{T}_3} + \mathbb{E}\|\widehat{g}-g_0\|_n^2 1_{\mathcal{T}_3^c} \qquad &\le O(\mu^2+\lambda^2 s_0) + O(1/\mu^2)\sqrt{\exp(-nc\rho^2)} \\ \le O(\mu^2+\lambda^2 s_0),&\end{aligned}$$

where the last step is true since $\sqrt{\exp(-n^{1/(2\alpha+1)})} = O(n^{-4/(2\gamma+1)})$.

Finally it follows from Lemma 4.1 of Nussbaum (1985) that $\mathbb{E}\int_0^1 |\widehat{f}(u) - f_0(u)|^2\, du = O(\mu^2+\lambda^2 s_0)$ and $\mathbb{E}\int_0^1 |\widehat{g}_0(u) - g_0(u)|^2\, du = O(\rho^2+\lambda^2 s_0)$.

□

## A.3. Results from empirical process theory

In this section, we include two theorems of empirical process theory, A.11 and A.12, which are Theorems 2.1 and 3.1 from van de Geer (2014), respectively.

Let $R_1^* = \sup_{f\in\mathcal{F}^*}\|f\|$, $K_1^* = \sup_{f\in\mathcal{F}^*}\|f\|_\infty$, and $R_2^* = \sup_{g\in\mathcal{G}^*}\|g\|$, $K_2^* = \sup_{g\in\mathcal{G}^*}\|g\|_\infty$.

**Theorem A.11.** *For all $t>0$, with probability at least $1-\exp(-t)$,*

$$\sup_{f\in\mathcal{F}^*}\Big|\|f\|_n^2 - \|f\|^2\Big|/C_1 \le \frac{2R_1 J_\infty(K_1^*, \mathcal{F}^*) + R_1^* K_1^* \sqrt{t}}{\sqrt{n}} + \frac{4J_\infty^2(K_1^*, \mathcal{F}^*) + K_1^{*2} t}{n}$$

*for some constant $C_1>0$.*

**Theorem A.12.** *Suppose that $R_1^*/R_2^* \le K_1^*/K_2^*$. For any $t \ge 4$ and $n$ such that*

$$\frac{2R_1^* \mathcal{J}_\infty(K_1^*, \mathcal{F}^*) + R_1^* K_1^* \sqrt{t}}{\sqrt{n}} + \frac{4\mathcal{J}_\infty^2(K_1^*, \mathcal{F}^*) + K_1^{*2} t}{n} \le \frac{R_1^{*2}}{C_1}$$

*and*

$$\frac{2R_2^* \mathcal{J}_\infty(K_2^*, \mathcal{G}^*) + R_2^* K_2^* \sqrt{t}}{\sqrt{n}} + \frac{4\mathcal{J}_\infty^2(K_2^*, \mathcal{G}^*) + K_2^{*2} t}{n} \le \frac{R_2^{*2}}{C_2},$$

*we have with probability at least $1-12\exp(-t)$,*

$$\frac{1}{8C_1}\sup_{f\in\mathcal{F}^*, g\in\mathcal{G}^*}\Big|(\mathbb{P}_n - \mathbb{P})fg\Big| \le \frac{R_1^* \mathcal{J}_\infty(K_2^*, \mathcal{G}^*) + R_2^* \mathcal{J}_\infty(R_1^* K_2^*/R_2^*, \mathcal{F}^*) + R_1^* K_2^* \sqrt{t}}{\sqrt{n}} + \frac{K_1^* K_2^* t}{n}.$$